\documentclass{article}
\usepackage[utf8]{inputenc}
\usepackage{amssymb}
\usepackage{amsfonts}
\usepackage{mathrsfs}
\usepackage{amsmath}
\usepackage{appendix}
\usepackage{graphicx}
\usepackage{dsfont}
\usepackage{float}
\usepackage{diagbox}

\usepackage{graphics}
\usepackage{csquotes}
\usepackage[all]{xy}
\usepackage[margin=1.4in,marginparwidth=1in,marginparsep=.05in]{geometry}
\PassOptionsToPackage{hyphens}{url}
\usepackage{tikz-cd}
\usepackage{braket}
\usepackage{xcolor}
\usepackage{nicematrix}
\usepackage{verbatim}
\usepackage{float}
\usepackage{mathrsfs}
\usepackage{extarrows}
\usepackage{tikz}
\usetikzlibrary{shapes,arrows,positioning}
\usepackage{mathrsfs}
\usepackage{mathtools}

\usepackage{subfigure}
\usepackage{hyperref}
\usepackage{wrapfig}
\usepackage{gensymb}
\usepackage{tabularx, booktabs, makecell, caption}

\usepackage{algorithm} 
\usepackage{algpseudocode} 
\usepackage{fancyhdr}

\begin{document}
\title{Efficient Spectral Element Method for the Euler Equations on Unbounded Domains}
\author{Yassine Tissaoui\footnote{Department of Mechanical and Industrial Engineering, New Jersey Institute of Technology, Newark NJ, U.S.A.. tissaoui@wisc.edu}, James F. Kelly\footnote{Space Science Division, U.S. Naval Research Lab, Washington DC, U.S.A.}, Simone Marras\footnote{Department of Mechanical and Industrial Engineering, New Jersey Institute of Technology, Newark NJ, U.S.A.. smarras@njit.edu}}
\date{\today}
\maketitle

\begin{abstract}
Mitigating the impact of waves leaving a numerical domain has been a persistent challenge in numerical modeling. Reducing wave reflection at the domain boundary is crucial for accurate simulations. Absorbing layers, while common, often incur significant computational costs. This paper introduces an efficient application of a Legendre-Laguerre basis for absorbing layers for two-dimensional non-linear compressible Euler equations. The method couples a spectral-element bounded domain with a semi-infinite region, employing a tensor product of Lagrange and scaled Laguerre basis functions. Semi-infinite elements are used in the absorbing layer with Rayleigh damping. In comparison to existing methods with similar absorbing layer extensions, this approach, a pioneering application to the Euler equations of compressible and stratified flows, demonstrates substantial computational savings. The study marks the first application of semi-infinite elements to mitigate wave reflection in the solution of the Euler equations, particularly in nonhydrostatic atmospheric modeling. A comprehensive set of tests demonstrates the method's versatility for general systems of conservation laws, with a focus on its effectiveness in damping vertically propagating mountain gravity waves, a benchmark for atmospheric models. Across all tests, the model  presented in this paper consistently exhibits notable performance improvements compared to a traditional Rayleigh damping approach.
\end{abstract}

\section{Introduction}

The challenge of preventing the reflection of waves as they exit a bounded domain has persisted for decades; among the very first articles describing this problem we mention, e.g., \cite{engquistMajda1977, engquistMajda1979}, with \cite{benacchioBonaventura2019, vismaraBenacchio2023} among the most recent. To resolve this issue there are two commonly used approaches. The first is to use radiative/non-reflecting boundary conditions which typically rely on a simplification of the full model dynamics to allow perturbations to pass through the boundary. We mention \cite{engquistMajda1977,engquistMajda1979,israeli1981approximation,neta2008application,giles1990} as some examples of this approach.  For example, Higdon boundary conditions \cite{higdon1986absorbing,higdon1994radiation,dea2011experimental} assume the outgoing waves have fixed wave speeds at the non-reflecting boundary, which poses a problem for nonlinear and dispersive waves present in the compressible Euler equations \cite{deaGiraldoNeta2009}.  Open boundary conditions for gravity waves \cite{klempDurran1983} perform well for idealized test cases, but are not sufficiently robust for real forecasts.  More recently, open BCs for high-altitude applications \cite{klemp2022constant,kelly2023physics} have been proposed to model the expansion and contraction of the thermosphere.  The main drawback of all of these approaches is that they approximate the model equations near the boundary.  As a result, many are only effective for linearized or idealized cases. 

A second and more robust approach is Rayleigh damping, which relaxes one or more of the state variables in an absorbing layer near the boundary.
This approach augments the computational domain by adding
additional levels/cells  in the absorbing layer 
that damps outgoing waves or disturbances toward a specified reference (see, e.g. \cite{durranKlemp1983,klemp1978numerical,klempDudhia2008,lavelle2008pretty}.) The effectiveness of these methods often depends on a tunable relaxation coefficient and thickness of the absorbing layer.  While many profiles have been employed, often the absorbing layer's thickness and parameters are determined empirically \cite{Modave2010}. Since the absorbing region must be discretized with additional cells or elements to appropriately damp outgoing waves, these additional cells/elements incur additional computational and communication cost for large, high-resolution, multidimensional problems. 

As a cost-effective alternative to stadard Rayleigh damping, Benacchio and Bonaventura  \cite{benacchioBonaventura2013} proposed a spectral collocation approach based on scaled Laguerre functions to prevent the reflection of one-dimensional shallow water waves. This method was subsequently extended by the same authors \cite{benacchioBonaventura2019} to the discontinuous Galerkin (DG) method. The governing equations are discretized with Lagrange polynomials in the finite region of the domain, while a scaled Laguerre basis is used in the semi-infinite region. The approach was further applied to solve the advection-diffusion equation on unbounded domains \cite{vismaraBenacchio2022}, and extended to address hyperbolic and parabolic problems in two dimensions \cite{vismaraBenacchio2023}.
For the continuous spectral element solution of the Stokes and Navier-Stokes equations governing incompressible flows, a coupled Legendre-Laguerre solution was proposed in \cite{zhuang2010coupled}, building on the earlier work presented in \cite{shen2000stable,shen2009some}. 

The current study uses the nodal continuous Galerkin (CG) combined with semi-infinite elements to solve the compressible nonlinear Euler equations for with gravity. This paper couples a bounded domain discretized using nodal continuous Galerkin spectral elements with tensor product bases, to a semi-infinite domain discretized using a basis that is the tensor product of a Lagrange polynomial basis on Legendre-Gauss-Lobatto (LGL) points in the normal direction  to the outgoing flow, and a Laguerre function basis on Laguerre-Gauss-Radau (LGR) points along the direction of the outgoing wave. Like \cite{zhuang2010coupled,benacchioBonaventura2019,vismaraBenacchio2022,vismaraBenacchio2023}, this approach neither seeks to find an analytical expression for the transmissive boundary conditions as done in e.g.\cite{deaGiraldoNeta2009,kelly2023physics}, nor does it attempt to be a selective damping approach that targets specific wave lengths for damping e.g.\cite{abarbanel1999well,navon2004perfectly,berenger1994perfectly,hesthaven1998analysis}. The use of the Laguerre function basis on semi-infinite elements allows an effective approximation of exponential decay and lends itself to the construction of efficient and effective Rayleigh damping layers. Moreover, the use of nodal continuous Galerkin elements with tensor product bases allows us to use computationally efficient algorithms while making the coupling between finite and semi-infinite elements seamless. 

To the authors' knowledge, this paper presents the first application of Laguerre-based absorbing layers to the two-dimensional solution of the Euler equations of compressible flows. Due to the importance of effective absorbing layers in atmospheric simulations, this paper studies their effectiveness to damp outgoing gravity waves although their applicability is not limited to this problem. As noted by \cite{vismaraBenacchio2023}, the Laguerre absorbing layer is many times cheaper than a classical Rayleigh damping approach and scales favorably with respect to the absorbing layer thickness.

The remainder of this paper is organized as follows: \S \ref{sec:problem} presents the governing equations. \S \ref{sec:num} presents the nodal continuous Galerkin spectral and semi-infinite element method, the algorithms used to construct the solution on semi-infinite elements, and other numerical techniques that are implemented in {\tt Jexpresso} \cite{jexpressoGithub}, a new open-source spectral element code written in performant Julia. Numerical results for a hierarchy of 1D and 2D tests are presented in \S \ref{sec:num_res} (All of the test cases presented in this paper are documented and ready to run by the user of the code), followed by discussion in \S \ref{sec:disc} and conclusions in \S \ref{sec:conc}. The pseudo-code to compute the element right-hand side on semi-infinite elements is given in the Appendix A while Appendix B provides a description of how to extend the current method to three dimensions. Appendix C shows that for mountain wave problems, no significant errors are introduced by the method when compared to a standard spectral element approach.

\section{Governing equations}
\label{sec:problem}
This paper considers a general system of PDEs, written in conservation form, on the semi-infinite domain $\Omega$ for $0 \leq t \leq T_{\mbox{end}}$:
\begin{equation}
\label{eq:CL}
\frac{\partial \textbf{q}}{\partial t} + \frac{\partial\textbf{F}(\textbf{q})}{\partial x} + \frac{\partial\textbf{G}(\textbf{q})}{\partial z}  = \textbf{S}(\textbf{q}) + \textbf{V}(\textbf{q}) +\textbf{R}(\textbf{q}),
\end{equation}
where the state vector $\textbf{q}$, flux vectors $\textbf{F}$ and $\textbf{G}$, source vector $\textbf{S}$, and diffusion vector $\textbf{V}$ are problem-dependent. A Rayleigh damping layer $\textbf{R}$ is used to damp outgoing waves, where required, is defined as 
\begin{equation}
    \mathbf{R}(\mathbf{q}) = -\gamma (z) \left( \mathbf{q} - \mathbf{q}_0 \right),
\end{equation} 
where $\gamma (z)$ is a damping coefficient to be defined when needed and $\mathbf{q}_0$ is the reference state. The following equations are be solved in this paper:
\begin{enumerate}
    \item 1D wave equation:
\begin{equation}\label{eq:Wave_eq}
{\bf q}=\begin{bmatrix}
u \\
v
\end{bmatrix},\quad{\bf F}=\begin{bmatrix}
v\\
u
\end{bmatrix},
\end{equation}
where the wave-speed is one.
    \item 1D linear shallow water equations:
 \begin{equation}\label{eq:Shallow}
{\bf q}=\begin{bmatrix}
h \\
u
\end{bmatrix},\quad{\bf F}=\begin{bmatrix}
Uh + Hu\\
gh + Uu
\end{bmatrix},
\end{equation} 
where $h$ is height, $u$ is velocity, $g$ is the acceleration of gravity, and $H$ and $U$ are the reference height and velocity, respectively.
\item 2D Helmholtz equation:
    \begin{equation}\label{eq:Helmholtz}
{\bf S}=\begin{bmatrix}
\alpha^2 u + f(x,z)
\end{bmatrix},\quad{\bf V}=\begin{bmatrix}
u_{xx} + u_{zz}
\end{bmatrix},
\end{equation}
where $\alpha$ is the wavenumber and $f(x,z)$ is a source function.
    \item 2D scalar advection-diffusion equation:
\begin{equation}\label{eq:2D_AdvDiff}
{\bf q}=\begin{bmatrix}
q\\
\end{bmatrix},\quad{\bf F}=\begin{bmatrix}
qu\\
\end{bmatrix},\quad{\bf G}=\begin{bmatrix}
qv\\
\end{bmatrix},\quad{\bf V}=\nu\begin{bmatrix}
q_{xx} + q_{zz}
\end{bmatrix}.
\end{equation}

    \item 2D compressible Euler equations with gravity:
\begin{equation} \label{eq:2D_Euler}
{\bf q}=\begin{bmatrix}
\rho \\
\rho u\\
\rho v\\
\rho \theta
\end{bmatrix},\quad{\bf F}=\begin{bmatrix}
\rho u\\
\rho u^2 + p\\
\rho u v\\
\rho u \theta
\end{bmatrix},\quad{\bf G}=\begin{bmatrix}
\rho v\\
\rho v u\\
\rho v^2 + p\\
\rho v \theta
\end{bmatrix},\quad{\bf S}=\begin{bmatrix}
0\\
0\\
-\rho g\\
0
\end{bmatrix},\quad{\bf V}=\begin{bmatrix}
0\\
\mu(u_{xx} + u_{zz})\\
\mu(v_{xx} + v_{zz})\\
\kappa(\theta_{xx} + \theta_{zz})
\end{bmatrix},
\end{equation}
where $\rho$ is density, $u$ is horizontal velocity, $v$ is vertical velocity, $\theta$ is potential temperature, $p$ is pressure, and $\mu$ and $\kappa$ are constant artificial viscosity coefficients.
 
\end{enumerate}

To preserve a numerically balanced hydrostatic state of the background atmosphere when solving the Euler equations, the hydrostatic balance terms in the vertical momentum equation are replaced by 
\[
\frac{\partial p'}{\partial z} \quad \rm{ and \quad } \rho'g, \]
where the primed quantities $\alpha' = \alpha - \alpha_0$ are the perturbations with respect to a hydrostatically balanced background state $\alpha_0$ (see e.g. \cite{marrasEtAl2015review} and citations therein).  In addition, the proposed semi-infinite element approach requires that each component of the state vector decays to zero as the vertical coordinate $z \rightarrow \infty$.
While density decays to zero as $z \rightarrow \infty$,  neither velocity nor potential temperature approach zero at high altitudes; hence, perturbation variables allows the proposed method to respect this decay condition while maintaining a physically meaningful reference state.

The equation of state of ideal gases is used to close the system.  Written in terms of potential temperature $\theta$, the pressure is 
\begin{equation} \label{eq:eos}
p = p_0\left(\frac{\rho R \theta}{p_0} \right)^{c_p/c_v},
\end{equation}
where $p_0=1000~$hPa, $c_p = 1005.0~$J/(kg K), $c_v = 718.0~$J/(kg K) and the specific gas constant $R = c_p - c_v$.

\section{Numerical methods}
\label{sec:num}

\begin{figure}
    \centering
    \includegraphics[width=0.3\textwidth]{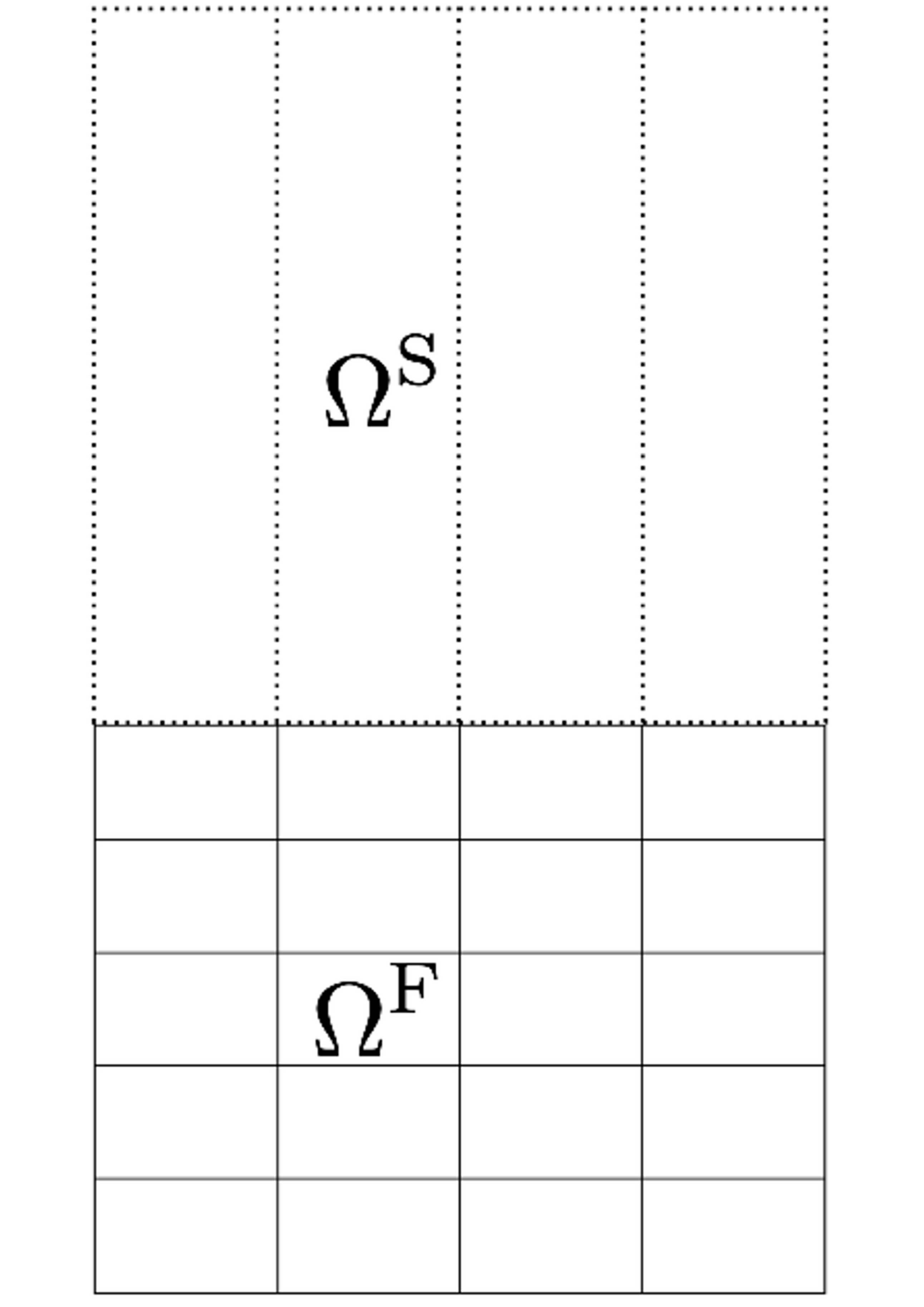}
    \caption{Example of a finite spectral element domain $\Omega^F$ with $N_F = 16$ connected to a semi-infinite element domain $\Omega^{S}$ with $N_S = 4$.}
    \label{fig:twodomains}
\end{figure}

The governing equations are discretized in space by means of high-order spectral elements with tensor-product Lagrange bases. In contrast to the modal discontinuous Galerkin approach of \cite{benacchioBonaventura2019}, we use nodal continuous Galerkin spectral elements. Given a computational domain $\Omega$, we first subdivide it into a finite domain $\Omega^F$ and a semi-infinite domain $\Omega^S$.  Each of these subdomains are then divided into
conforming quadrilaterals that are either finite spectral elements $\Omega^{F}_{e}$, or semi-infinite elements $\Omega^{S}_{e}$ to form the discrete domain $\Omega_h$ as

\begin{equation}
    \Omega \approx \Omega^h = \bigcup_{e=1}^{N_e} \Omega_e = \bigcup_{e=1}^{N_F} \Omega^{F}_e ~\cup~ \bigcup_{e=1}^{N_S} \Omega^{S}_e
\end{equation}
where $N_F$ and $N_{S}$ are the number of finite spectral elements and semi-infinite elements, respectively, and the total number of elements is $N_e = N_F + N_S$.  This configuration is illustrated in Figure \ref{fig:twodomains}.

In time, the governing equations are discretized by means of an explicit 3-stage, $3^{rd}$-order Runge-Kutta approximation. {\tt Jexpresso} \cite{jexpressoGithub}, the open-source spectral element code utilized for this study, relies on the {\tt DifferentialEquation.jl} package \cite{rackauckas2017differentialequations} for time integration.

\subsection{Basis functions on the finite domain}
The spectral element solution of partial differential equations relies on the interpolation of the unknowns and their derivatives. Interpolation is achieved by means of Lagrange polynomials $h(\xi)$ of degree $N$. 
For a given set of interpolating nodes $\xi_i$ ($i=1,\dots, N+1$), the $i^{th}$ polynomial is the orthogonal function 
\begin{equation}
\label{lagrangeDef}
h_i(\xi) = \prod_{j=0,j\neq i}^{N}\frac{\xi - \xi_j}{\xi_i - \xi_j}.
\end{equation}
By definition, $h_i$ evaluated at point $\xi_j$ simplifies to

\begin{equation}
    h_i(\xi_j) = \delta_{ij} =
    \begin{cases}
      1 & \text{if $i=j$}\\
      0 & \text{otherwise,}
    \end{cases}
\end{equation}
where $\delta_{ij}$ is the Kronecker delta and $\xi_j$ are the coordinates of the $j^{th}$ interpolating point.

The Legendre-Gauss-Lobatto (LGL) points are used both as interpolation and integration nodes (inexact integration).  The LGL points are the zeros of the $N$-th degree Lobatto polynomial, given by
\begin{equation}
    (1-\xi^2)P^{'}_N(\xi)=0,
\end{equation}
where $P_N(x)$ is the Legendre polynomial of order $N$. These polynomials and their first derivative are computed recursively via
\begin{subequations}
\label{legendreEq}
\begin{equation}
P_0(\xi) = 1
\end{equation}
\begin{equation}
P_1(\xi) =\xi
\end{equation}
\begin{equation}
P_{k}(\xi) = \frac{2k - 1}{k}\xi P_{k-1}(\xi) - \frac{k-1}{k}P_{k-2}(\xi),~~~~\forall k\geq 2\\
\end{equation}
\begin{equation}
P'_k(\xi) =  (2k -1)P_{k-1}(\xi) + P'_{k-2}(\xi).
\end{equation}
\end{subequations}
See, e.g. \cite{giraldoBOOK, karniadakisSherwinBook2005}.  Finally, the corresponding quadrature weights are given by
\begin{equation}
    \omega(\xi_i)=\frac{2}{N(N+1)}\left[\frac{1}{P_N(\xi_i)}\right]^2.
\end{equation}
%

For every spectral element, we seek an approximation $f^h$ of the variable $f$ by means of the expansion
\begin{equation}\label{eq:f_h}
    f^h(\boldsymbol{\xi},t) = \sum_{l=1}^{(N+1)^2} \psi_l(\boldsymbol{\xi})\hat{f_l}(t),
\end{equation}
where $\boldsymbol{\xi}=(\xi, \eta)$, $\hat{f_l}(t)$ are the expansion coefficients, and $\psi_l(\boldsymbol{\xi})$ are two-dimensional nodal basis functions defined as the tensor product of the 1D Lagrange polynomials
\begin{equation}\label{eq:tensor_product}
    \psi_l(\mathbf{x}) = h_i[\xi(\mathbf{x})] \otimes h_j[\eta(\mathbf{x})], \hspace{5pt} \hspace{3pt} l=i+1 + j(N+1),
\end{equation}
where $\mathbf{x}=(x,z)$. 

As such, we can approximate an integral over an element $\Omega_e$ using an inexact Gauss quadrature rule of accuracy $\mathcal{O}(2N-1)$ as follows:
\begin{equation}\label{eq:Gauss_quadrature}
    \int_{\Omega_e}f(\mathbf{x})d\mathbf{x} = \int_{\Omega_{ref}}f(\boldsymbol{\xi})|\mathbf{J}(\boldsymbol{\xi})|d\boldsymbol{\xi}\approx \sum_{i,j=1}^{N+1}\omega(\xi_i) \omega(\eta_j) f(\xi_i,\eta_j)|\mathbf{J}(\xi_i,\eta_j)|,
\end{equation}
where $|\mathbf{J}|$ is the determinant of the Jacobian matrix.
To keep the inexact quadrature error small, polynomials of order 3 or higher should be used.

\subsection{Basis functions on the semi-infinite domain}
The semi infinite domain uses both a different set of basis functions and integration points than its finite domain counterpart. First, define the Laguerre polynomials using their three term recurrence relation \cite[Eq. (3.3)]{shen2009some}:
\begin{subequations}
\label{laguerre}
\begin{equation}
    L_0(\xi) = 1 
\end{equation}
\begin{equation}
    L_1(\xi) = 1 - x 
\end{equation}
\begin{equation}
    L_k(\xi) = \frac{2k -1 -\xi}{k}L_{k-1}(\xi) - \frac{k-1}{k}L_{k-2}(\xi), \forall k \geq 2.
\end{equation}
\end{subequations}
Additionally the first derivative of the $k$-th Laguerre polynomial is \cite[Eq. (3.5)]{shen2009some}:
\begin{equation}
    L'_k(\xi) = -\sum_{n=0}^{k-1}L_{n}(\xi).
\end{equation}

The Laguerre polynomials are orthogonal on the semi-infinite interval $[0,\infty)$ with respect to an exponentially decaying weight as follows:

\begin{equation}\label{eq:Laguerre_ortho}
    \int_0^{\infty} L_i(\xi)L_j(\xi)e^{-\xi} d\xi = \delta_{ij} ,\textrm{$\forall i,j \geq 0$}.
\end{equation}

The Laguerre-Gauss-Radau (LGR) points are the roots of $\xi L'_{N+1}(\xi)$ for a fixed integer N. The LGR points $\{\xi_j\}_{j=0,N}$ will be used to construct a nodal spectral element on 
the semi-infinite elements.  We compute the LGR points using the Eigenvalue Method \cite{golub1969calculation,shen2009some}, which forms a tridiagonal matrix using the coefficients in the three-term recurrence relationship \eqref{laguerre} and solves an eigenvalue problem.  This method is stable and robust at very high-order and we have tested it for orders up to 60. 
Following \cite{shen2000stable}, the scaling factor $\lambda$ adjusts the LGR nodes, such that the physical nodes on the semi-infinite element are $\{x\}^{N+1}_{i=1} = \lambda \{\xi\}^{N+1}_{i=1}$. This scaling factor allows us to adjust the effective length of the semi-infinite element for a given problem.

\begin{figure}[h!]
    \centering
    \includegraphics[width=0.9\textwidth]{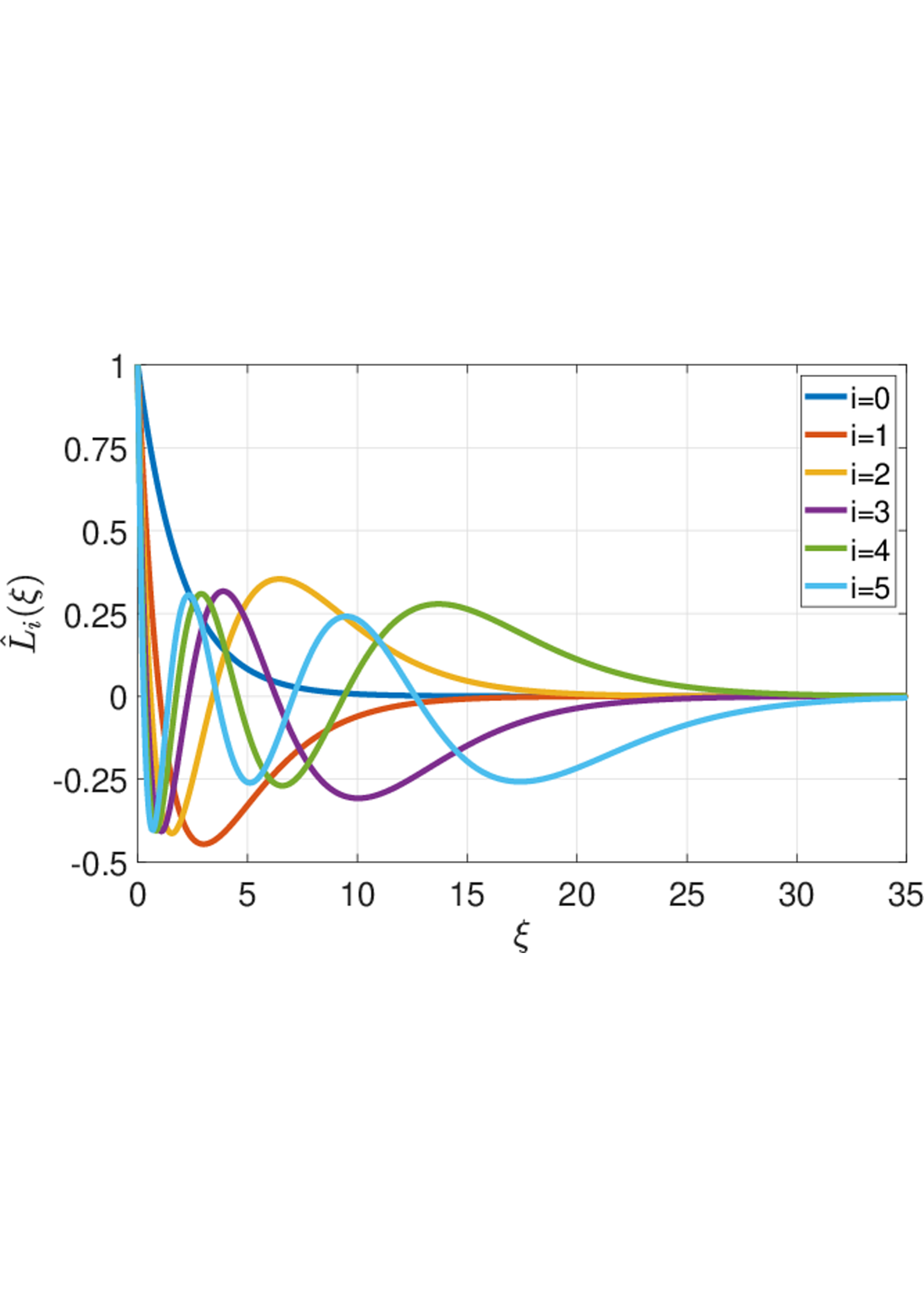}
    \caption{First six scaled Laguerre functions (SLFs) specified by \eqref{eq:laguerre_function} with scaling factor $\lambda=1$.}
    \label{fig:Laguerre_functions}
\end{figure}

Next we construct the Lagrange-Laguerre interpolating polynomials associated with the LGR points $\{\xi_j\}_{j=0,\dots,N}$ points following \cite{shen2000stable}:
\begin{equation}
    h_j^{Lag}(\xi) = -\frac{\xi L^{'}_{N+1}(\xi)}{(N+1)L_{N+1}(\xi_j)(\xi-\xi_j)}.
\end{equation}
We can then write their derivatives as follows:
\begin{equation}
    h_j^{'Lag}(\xi_i) = \begin{cases}
        \frac{L_{N+1}(\xi_i)}{L_{N+1}(\xi_j)(\xi_i-\xi_j)}  & \textrm{if $i \neq j$} \\
        \frac{1}{2} & \textrm{if $i = j \neq 0 $} \\
        \frac{N}{2} & \textrm{if $i = j = 0$.} 
        \end{cases}
\end{equation}
We now introduce the scaled Laguerre function (SLF) \cite{shen2000stable,shen2009some,benacchioBonaventura2013} 
\begin{equation}\label{eq:laguerre_function}
    \hat{L}_i(\xi) = e^{-\frac{\xi}{2\lambda}} L_i \left(\frac{\xi}{\lambda} \right),
\end{equation}
where $\lambda$ is a scaling factor and represents a characteristic length. Note that this notation is equivalent to the notation in \cite{benacchioBonaventura2013} for $\lambda = \beta^{-1}$. 
Applying \eqref{eq:Laguerre_ortho} yields
\begin{equation}
     \int_0^{\infty} \hat{L}_i(\xi)\hat{L}_j(\xi) d\xi = \lambda \delta_{ij}, \forall i,j \geq 0,
\end{equation}
indicating that the SLFs form an orthogonal basis on $L^2(\mathbb{R}^+)$. Each SLF decays  exponentially as $\xi \rightarrow \infty$ for any $\lambda>0$.  This property is illustrated in figure \ref{fig:Laguerre_functions}, which shows the first six Laguerre functions $\hat{L}_i(\xi)$ for $\lambda=1$.
Thanks to this property, the SLFs \eqref{eq:laguerre_function} are ideal for approximating functions in an absorbing layer; the layer damps any outgoing perturbations by enforcing the exponential decay property.

\emph{\textbf{Remark:}} We limit the use of damping terms to $\Omega^S$ and essentially overlap the sponge layer with the semi-infinite elements. In the results, we will describe the damping coefficients used for each test and we will show that minimal reflection can be obtained while relying on this approach.

We can now construct the Lagrange-Laguerre interpolant associated with LGR points $\{\hat{h}^{Lag}_j\}_{j=0}^N$ such that $\hat{h}^{Lag}_j(\xi_i)=\delta_{ij}$, and
\begin{equation}
    \hat{h}^{Lag}_j(\xi) = \frac{\exp(-\xi/2)}{\exp(-\xi_j/2)}h^{Lag}_j(\xi).
\end{equation}
Their derivatives are:
\begin{equation}
    \hat{h}_j^{'Lag}(\xi_i) = \begin{cases}
        \frac{\hat{L}_{N+1}(\xi_i)}{\hat{L}_{N+1}(\xi_j)(\xi_i-\xi_j)}  & \textrm{if $i \neq j$} \\
        0 & \textrm{if $i = j \neq 0 $} \\
        -\frac{N+1}{2} & \textrm{if $i = j = 0$.} 
        \end{cases}
\end{equation}
The quadrature weights $\{\hat{\omega}(\xi_i)\}^N_{i=0}$ associated with the LGR points are defined as
\begin{equation}
    {\hat\omega(\xi_i)} = \frac{\exp(\xi_i)}{(N+1)[L_N(\xi_i)]^2}.
\end{equation}

Let us now consider the reference 2D semi-infinite element $\Omega^S_{ref}$ such that it is only semi-infinite in the direction of the outgoing waves. In its finite direction, the LGL nodes and their associated bases are used, and in its semi-infinite direction the LGR nodes and their associated SLF bases are used. As such, for an integration point $\mathbf{\xi} \in \Omega^S_{ref}$, we can write $\mathbf{\xi} = (\xi_{LGL},\eta_{LGR})$.
The nodal basis functions $\psi^S_l$ on the semi-infinite element are the tensor product of the 1D Lagrange polynomials associated with the LGL nodes and the 1D Lagrange-Laguerre interpolating functions associated with the LGR nodes giving us the following:
\begin{equation}\label{eq:basis}
    \psi^S_l(\mathbf{x})_{\Omega^S_{ref}} = h_i[\xi(\mathbf{x})] \otimes \hat{h}_j[\eta(\mathbf{x})], \hspace{5pt} \hspace{3pt} l=i+1 + j(N_{LGL}),
\end{equation}
where $i \in \{1,\dots,N_{LGL}\}$, $j \in \{1,\dots,N_{LGR}\}$, $N_{LGL}$ is the number of LGL nodes and $N_{LGR}$ is the number of LGR nodes. This makes it so that equation \eqref{eq:Gauss_quadrature} remains valid on semi-infinite elements provided the appropriate substitutions of nodes and weights is performed. 

\subsection{Constructing element matrices}
In this section we present the reader with a template for how an element matrix is constructed for elements of the finite and semi-infinite domains. In what follows we will discuss the construction of the mass matrix. We refer the reader to \cite{giraldoBOOK} for constructions of the differentiation or laplacian matrices on the finite domain. The extension of these constructions to the semi-infinite domain is done similarly to the construction of the mass matrix. Let us first define every component of the mass matrix on a given spectral or semi-infinite element $\Omega_e \in \Omega^h = \Omega^F \cup \Omega^S$:

\begin{equation}
\mathbf{M}_{ij}^e=\int_{\Omega_e}\psi_i(\mathbf{x}) \psi_j(\mathbf{x}) d \mathbf{x} = \int_{\Omega_{ref}} \psi_i(\boldsymbol{\xi})\psi_j(\boldsymbol{\xi})| \mathbf{J}(\boldsymbol{\xi})| d \boldsymbol{\xi} \textrm{  $\forall i,j = 1,\dots,N_{\xi}N_{\eta}$},
\end{equation}
where $\psi_i$ is defined using a generalized form of \eqref{eq:tensor_product} as follows:

\begin{equation}\label{eq:Mass_Integral}
    \psi_l(\mathbf{x}) = h_i[\xi(\mathbf{x})] \otimes \overline{h}_j[\eta(\mathbf{x})], \hspace{5pt} \hspace{3pt} l=i + (j-1)N_{LGL},
\end{equation}
and where,

\begin{equation}
    \Omega_{ref},N_{\xi}, N_{\eta}, \overline{h}_j = \begin{cases}
        \Omega^F_{ref}, N_{LGL}, N_{LGL}, h_j & \textrm{if $\Omega_e \in \Omega^F$} \\
        \Omega^S_{ref}, N_{LGL}, N_{LGR}, \hat{h}^{Lag}_j & \textrm{if $\Omega_e \in \Omega^S$}
    \end{cases}.
\end{equation}
By approximating \eqref{eq:Mass_Integral} using an inexact quadrature rule (i.e. the quadrature and interpolation points coincide) we obtain the following:
\begin{align}
    \mathbf{M}^e_{ij} &=\sum_{k=1}^{N_{\xi}} \sum_{m=1}^{N_{\eta}}\overline{\omega}(\xi_k,\eta_m)
    \psi_i(\xi_k,\eta_m) \psi_j(\xi_k,\eta_m)
    |\mathbf{J}(\xi_k,\eta_m)|,\\
    \overline{\omega}(\xi_k,\eta_m) &= \omega_{\xi}(\xi_k)\omega_{\eta}(\eta_m),\\
    \omega_{\xi},\omega_{\eta} &= \begin{cases}
        \omega, \omega & \textrm{ if $\Omega_e \in \Omega^F$}\\
        \omega, \hat{\omega} & \textrm{ if $\Omega_e \in \Omega^S$}
    \end{cases},\\
    N_{\xi},N_{\eta}
    &= \begin{cases}
        N_{LGL},N_{LGL} & \textrm{ if $\Omega_e \in \Omega^F$}\\
        N_{LGL},N_{LGR} & \textrm{ if $\Omega_e \in \Omega^S$}\\
    \end{cases}.
\end{align}

Next, we present the pseudo-code for constructing the mass matrix of an element in the semi-infinite domain using inexact integration. Let us define $\mathbf{M}^{Lag}$ as the mass matrix of an element of the semi-infinite domain, $\eta^{LGR}$ as the LGR nodes, and $\xi$ as the LGL nodes:
\begin{algorithm}[h!]
    \caption{Construction of the mass matrix of a semi-infinite element}
    \begin{algorithmic}
    \State $\mathbf{M}^{e,Lag} = {\rm zeros}(N_{LGL}N_{LGR},N_{LGL}N_{LGR})$ 
        \For {$l=1,N_{LGR}$}
        \For {$k=1,N_{LGL}$}
        \State $I = k + (l-1)(N_{LGL})$
        \State $\overline{\omega}=\omega(\xi_k)\hat{\omega}(\eta^{LGR}_l)$
        \State $\mathbf{x}=(\xi_k,\eta^{LGR}_l)$
        \For {$j=1,N_{LGR}$}
        \For {$i=1,N_{LGL}$}
        \State $J = i + (j-1)(N_{LGL})$
        \State $\mathbf{M}^{e,Lag}_{IJ}= \mathbf{M}^{e,Lag}_{IJ}+\overline{\omega}\psi_I(\mathbf{x})\psi_J(\mathbf{x})|\mathbf{J}(\xi_k,\eta^{LGR}_l)|$
        \EndFor
        \EndFor
        \EndFor
        \EndFor
    \end{algorithmic}
\end{algorithm}

\subsection{Direct stiffness summation}
\label{sec:dss}
In order to couple the element local Galerkin expansion given by \eqref{eq:f_h} between adjacent elements, we need to 
construct a direct stiffness summation (DSS) operator.  The DSS operator
enforces the continuity of the global solution by averaging the state variable on nodes shared by multiple elements.  As shown below, this DSS operator couples the finite domain $\Omega^F$ and the semi-infinite domain $\Omega^S$ illustrated in Figure \ref{fig:twodomains} in a straight-forward manner.

First, we must define mappings from local elements to global nodes.
Let $I = H^F(e,i)$ be the map from the local element-wise node $i$ on the $e$-th finite element $\Omega_e^F$ and let $I = H^S(e,i)$ be the corresponding map from the 
$e$-th semi-infinite element $\Omega_e^S$.  For $H^F$, $i$ runs from 1 to $N_{LGL}^2$, while for $H^S$, $i$ runs from 1 to $N_{LGL}N_{LGR}$.  These mappings contain the connectivity information in the finite and semi-infinite grids, respectively.

We can now illustrate how the global problem is assembled by considering the mass matrix.  Let $\mathbf{M}_i^{e}$ be a local diagonal mass matrix corresponding to element $e$ (either finite or semi-infinite).  We construct the global mass matrix via a DSS operator 
\begin{equation}
\mathbf{M}_{IJ} = \bigwedge_{e=1}^{N_e} \mathbf{M}^{e}_{ij} ,
\label{matvec2}
\end{equation}
where the DSS operator $\bigwedge_{e=1}^{N_e}$ consists of a local-global mapping and appropriate summation.  For additional details, see Section 5.8 
in \cite{giraldoBOOK} or \cite{kellyGiraldo2012} for the parallel MPI implementation.  This DSS operator may be decomposed into two independent DSS operators, the first over the collection of $N^f$ finite elements, and the second over the 
collection of $N^s$ semi-infinite elements via
\begin{equation}
\mathbf{M}_{IJ} = \left[ \bigwedge_{e_f=1}^{N_f} \mathbf{M}^{e_f}_{ij} \right] \bigwedge \left[ \bigwedge_{e_s=1}^{N_s} \mathbf{M}^{e_s}_{ij} \right] .
\label{eq:matvec3}
\end{equation}
Hence, the DSS operation consists of three stages: 1) perform a DSS over the finite domain $\Omega^F$ using the local to global mapping $H^F$, 2) perform a DSS over the semi-infinite domain $\Omega^S$ using the mapping $H^S$, and finally 3) DSS the nodes shared by both the finite and semi-infinite domains $\Omega^F \cap \Omega^S$ using both $H^F$ and $H^S$.  Since the only coupling between the finite and semi-infinite grids is via this final DSS ope rator, the proposed semi-infinite approach may be retrofitted to an existing spectral element solver with only minor modifications.

\paragraph{Low pass spectral filter}
To control unresolved grid-scale noise and aliasing, low-pass spectral filters are typically used with the spectral element method (SEM) [Sec. 18.3]\cite{giraldoBOOK}.  Spectral filters employ a three-step process to damp/remove unphysical high-frequency components: 1) the element-local nodal solution is transformed into modal space, 2) a low-pass filter is applied in modal space, 3) the filtered modal representation is inverse transformed to nodal space.  The nodal to modal transform requires defining a set of modal basis functions.  As discussed in Section 3d of \cite{giraldo2004scalable}, the low-pass filter should not violate the continuity requirement of SEM.  One way to enforce this is by choosing modal functions such that most of the modes are zero at the element boundary.  For a 1D element using LGL points, an appropriate choice is $\phi_k(\xi) = P_k (\xi)$ for $k= 0$ or 1 and $\phi_k (\xi) = P_k(\xi) - P_{k-2} (\xi)$ for $k \geq 2$.  Since $\phi_k(\pm 1) = 0$ for $k \geq 2$, these higher-order modes do not effect the boundary of the element, and hence may be damped by an appropriate filter function.  The resulting transform (or Legendre) matrix is given by Eq. (31) in \cite{giraldo2004scalable}.  

For a 1D semi-infinite element using LGR points, a similar choice of modal functions is $\phi_0 (\xi) = e^{-\xi/2}$ and $\phi_k (\xi) e^{-\xi/2} \left[ L_k (\xi) - L_{k-1} (\xi) \right]$.  Since $L_k (0) = 1$ for all $k \geq 0$, we have $\phi_k(0) = 0$ for $k \geq 1$.  The Legendre matrix is then constructed by evaluating $\phi_k (\xi_j)$, where $\xi_j$ are the LGR points.  Since we employ a tensor product of LGL and LGR points within semi-infinite element, we construct a corresponding 2D tensor product of modal functions in order to transform the nodal representation.  A Boyd-Vandeven filter \cite{boyd1996erfc} is applied to the modal representation, and then the solution is inverse transformed to nodal space.

\section{Numerical Results}\label{sec:num_res}

We assess the efficacy and efficiency of the semi-infinite element approach with a series of one and two-dimensional tests. In 1D, we conduct two tests. The first test evaluates the method’s ability to damp outgoing waves using the wave equation, while the second test is a wave train for the linearized shallow water equations that demonstrates how the absorbing layer handles a series of incoming waves while maintaining the finite domain solution quality.
In 2D, we initially test the proposed method on the advection-diffusion and Helmholtz equations. For atmospheric simulations, we employ three standard idealized tests: the mountain-triggered hydrostatic gravity waves over a single mountain, and the non-hydrostatic wave problem triggered by a Sc\"ar ridge. These tests --with increasing physical complexity-- allow us to verify that the proposed approach is effective as an absorbing layer for the compressible Euler equations and showcase its potential for use in atmospheric simulations.

\subsection{1D wave equation}

We solve the 1D wave equation as described by \eqref{eq:CL} for the unknown and flux vectors \eqref{eq:Wave_eq}. The space considered is the finite domain $\Omega^F = [-2.5,2.5]~$m which is discretized using 50 spectral elements of order 6 leading to an effective resolution $\Delta x = 0.016~$m. One semi-infinite element is added on each end of the finite domain such that the left region is $\Omega^S_L = (-\infty ,-2.5]~$m, whereas the right is $\Omega^S_R = [2.5,\infty)~$m. We initialize $u$ and $v$ as
\begin{align}
    u(0,x) &=2^{-\frac{(x-x_c)^2}{\sigma^2}},\\
    v(0,x) &= 0,
\end{align}
where $x_c = 0$ and $\sigma = 0.15$ m.
While the scaled Laguerre functions (SLFs) present an excellent tool for accurately approximating exponentially decaying functions, they are not capable of acting as an absorbing layer on their own. To this end, a sponge is added to each Laguerre semi-infinite element to damp incoming waves while allowing the SLF basis to accurately approximate this damping. This allows a single semi-infinite element with a sponge to act as the absorbing layer on each side of the domain. We define $X_L$ and $X_R$ to be the leftmost and rightmost nodes of the domain respectively, such that $X_L =-11.63$ m and $X_R = 11.63$ m for semi infinite elements of order 50 and a scaling factor  $\lambda=0.05~$m. This test is run for 9 seconds of simulated time with a time step $\Delta t = 0.001~$s. The Rayleigh damping on each side is controlled by $\gamma(x)$, adapted from \cite{benacchioBonaventura2013}, defined for all $x \in \Omega^S_L \cup \Omega^S_R$ as follows:
\begin{equation}\label{eq:1D_damping}
    \gamma(x) = \frac{\Delta \gamma}{1+\exp\left(\frac{{\rm sgn}(L_0-X_0)(\alpha(L_0-X_0)-x+X_0)}{\zeta}\right)},
\end{equation}
where ${\rm sgn}$ is the sign function, $X_0$ is the coordinate of the first node of the semi-infinite element (in this case $X_0=-2.5$ m for the left element and $X_0=2.5$ m for the right element). $\Delta \gamma=2 \mbox{~s}^{-1}$, $L_0$ is the end point of each element, such that $L_0=X_L$ for the left element and $L_0=X_R$ for the right element, $\alpha=0.3$, and $\zeta=\max(|X_L|,|X_R|)/18$.

\begin{figure}[h!]
    \centering
    \includegraphics[width=0.99\textwidth]{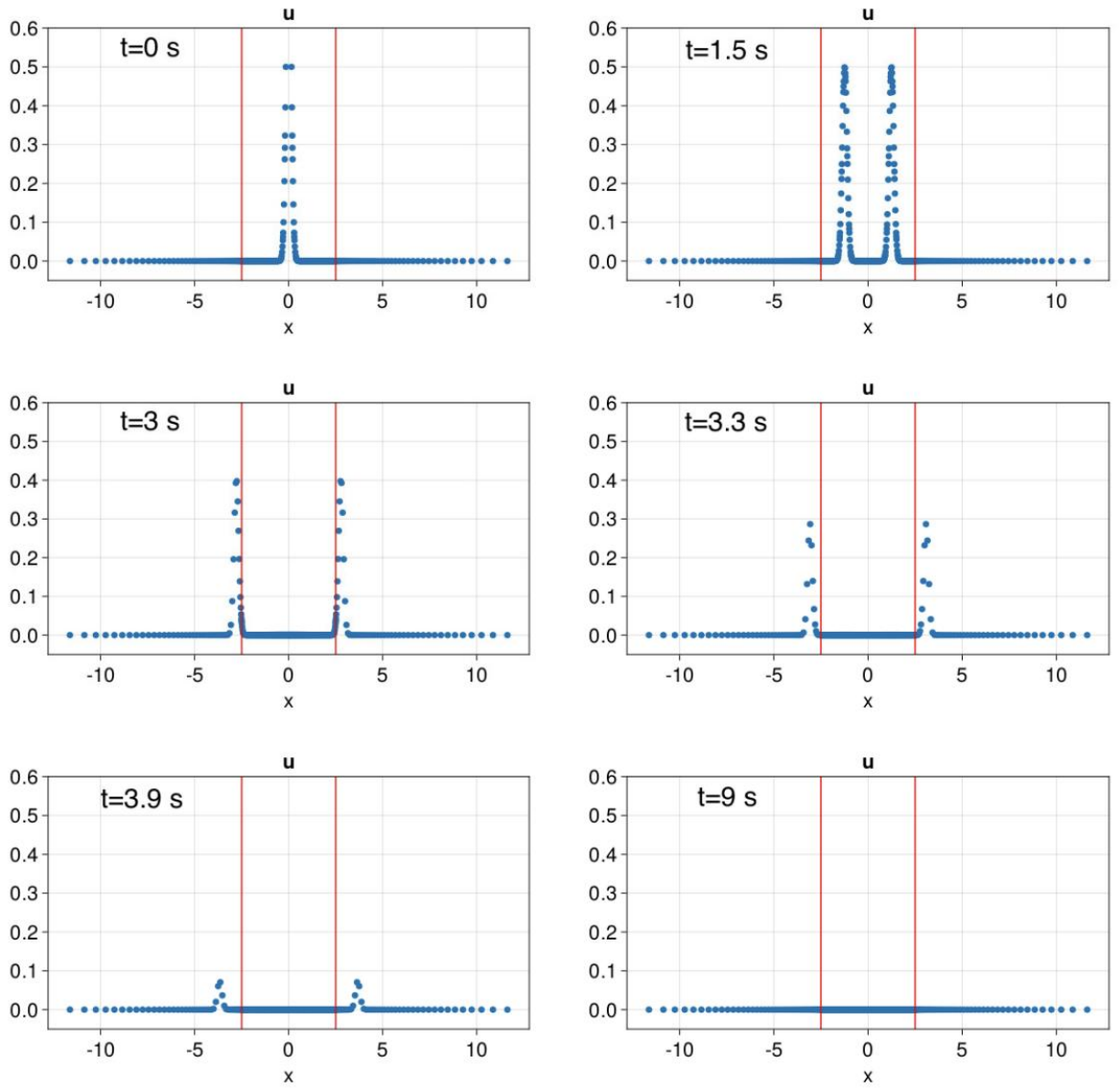}
    \caption{1D wave equation: Time evolution of $u$: Top-left: $t=0~$s, top-right: $t=1.5~$s, center-left: $t=3~$s, center-right: $t=3.3~$s, bottom-left: $t=3.9~$s, bottom-right: $t=9~$s. The initial perturbation splits into two waves, one left-going wave and one right-going wave, both reach the edges of the finite domain (illustrated by the red vertical lines) and are subsequently damped as they move through the semi-infinite elements.}
    \label{fig:Wave_u}
\end{figure}

\begin{figure}[h!]
    \centering
    \includegraphics[width=0.99\textwidth]{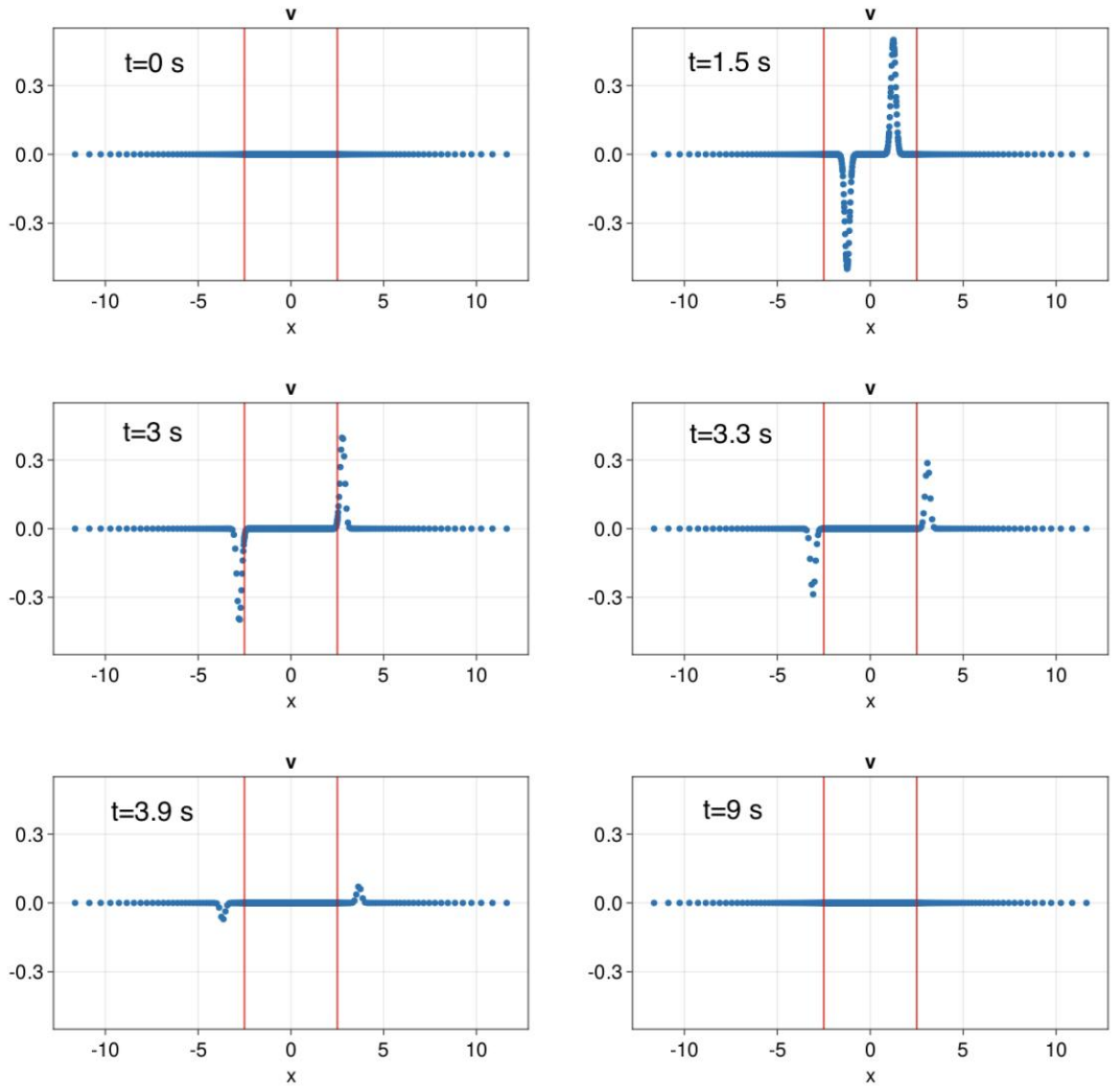}
    \caption{1D wave equation: Same as figure \ref{fig:Wave_u} but for $v$}
    \label{fig:Wave_v}
\end{figure}

As shown in figures \ref{fig:Wave_u} and \ref{fig:Wave_v}, the initial wave splits into two, one left-going wave with a negative velocity $v < 0$ and one right-going wave with a positive velocity $v > 0$, both eventually reach the edges of the finite domain and are then damped at the same rate in the semi-infinite absorbing layer.

Table \ref{tab:1D_Wave_Time} compares the cost of using the Laguerre semi-infinite elements with the cost of a standard Rayleigh damping using an extended finite domain.  These results are presented in the form of a non-dimensional time per time step $T^*$, which is simply the ratio of the times to solution of a given test case and the time to solution of the first test case in the table. We also present $T^*_{Finite}\%$, and $T^*_{Laguerre}\%$, the percentages of $T^*$ spent on the finite and semi-infinite respective right hand sides. For the tests conducted with no semi-infinite elements, the finite domain is extended to be $\Omega^F = [X_L,X_R]~$m. The spaces $[X_L,-2.5]~$m and $[2.5,X_R]~$m which would each be occupied by a semi-infinite element are instead discretized using enough finite domain elements of the same order (in this case the order is maintained at 6) to maintain the same resolution $\Delta x = 0.016~$m in the entire finite domain. In the table $N_x$ tracks the number of elements in $\Omega^F$ for each simulation.

We can see from the table that extending the finite domain while maintaining its resolution is considerably more time consuming than using semi-infinite elements. It is twice as costly compared to using order 20 Laguerre elements (for the same total domain extent), and nearly three times as costly when trying to mimic the order 50 Laguerre elements. Furthermore, even while using up to order 50 semi-infinite elements, the cost of the Laguerre right hand side remains lower than the finite domain right hand side while presenting a much higher order approximation. Additionally, this relative cost decreases further the more elements are present in $\Omega^F$ and we demonstrate this with the next test.

\begin{table}[h!]
    \centering
    \begin{tabular}{|c|c|c|c|c|c|c|}
    \hline
         Absorbing layer type& $T^*$ & $X_{L}$ & $X_{R} $&$T^*_{Finite} \%$& $T^*_{Laguerre} \%$ & $N_x$\\
    \hline
    \hline
         Semi-infinite element of order 20 & 1& -5.91& 5.91 &78.74& 21.26 & 50\\
         Semi-infinite element of order 25 & 1.07& -6.85& 6.85 &73.57& 26.43 & "\\
         Semi-infinite element of order 30 & 1.15& -7.8 & 7.8 &68.09& 31.91& "\\
         Semi-infinite element of order 35 & 1.24& -8.75 & 8.75& 63.41& 36.59 & "\\
         Semi-infinite element of order 40 & 1.33& -9.7& 9.7 &58.87& 41.13& "\\
         Semi-infinite element of order 50 & 1.53& -11.63& 11.63 &51.18& 48.82& "\\
         \hline
         Extended finite domain of order 6 & 2.05& -5.91& 5.91 &100& N/A & 118\\
         " & 2.36& -6.85& 6.85 & "& "& 137\\
         " & 2.74& -7.8& 7.8& "& "& 156\\
         " & 3.14& -8.75 & 8.75 &"& "& 175\\
         " & 3.45& -9.7& 9.7 &"& "& 194\\
         " & 4.2& -11.63 & 11.63 & "& "& 233\\
        \hline
    \end{tabular}
    \caption{Timings of the 1D wave equation simulations with and without Laguerre semi-infinite elements in the absorbing layer. Simulations using the extended finite domain have the same physical domain size as simulations using semi-infinite elements, but the space that would be otherwise occupied by these elements is instead occupied by elements of the finite domain retaining the same resolution $\Delta x = 0.016~$m and element order $N=6$.}
    \label{tab:1D_Wave_Time}
\end{table}

\subsection{1D wave train for linearized shallow water equations}
The wave train is found as the solution to \eqref{eq:CL} for the solution and flux vectors \eqref{eq:Shallow}.
This system describes the evolution of small perturbations of the free surface of a non-rotating fluid of height perturbation $h$ and velocity $u$ under constant gravity. 

We define the finite domain as $\Omega^F =[0,5000]~$m and initialize the fluid to be at rest such that $h(0,x) = 0~$m and $u(0,x)=0~{\rm ms^{-1}}$. The reference height and velocity are, respectively, $H=10~$m and $U=0~{\rm ms^{-1}}$. The flow is forced through a time dependent Dirichlet boundary condition at $x=0$ such that $u(t,0) = A\sin(2\pi k t/T)$,
where $T=5,000~$s is the final time of the simulation, $A=0.025~{\rm ms^{-1}}$ is the amplitude of the perturbation, and $k=30$. The finite domain is discretized using 300 spectral elements of order 4, yielding an average spatial resolution of $\Delta x= 4.16~$m. A Laguerre semi-infinite element is defined as $\Omega^S = [5000,\infty)~$m and the same damping coefficient defined in \eqref{eq:1D_damping} is used $ \forall x \in \Omega^S$, such that $L_0=X_{end}$, and $\zeta = L_0/18$. 
A Laguerre semi-infinite element of order 50 is used in the damping layer which, with a scaling factor $\lambda =100~$m, results in the end point $X_{end}=23,262~$m.

\begin{figure}[h!]
	\centering
    \includegraphics[width=0.99\textwidth]{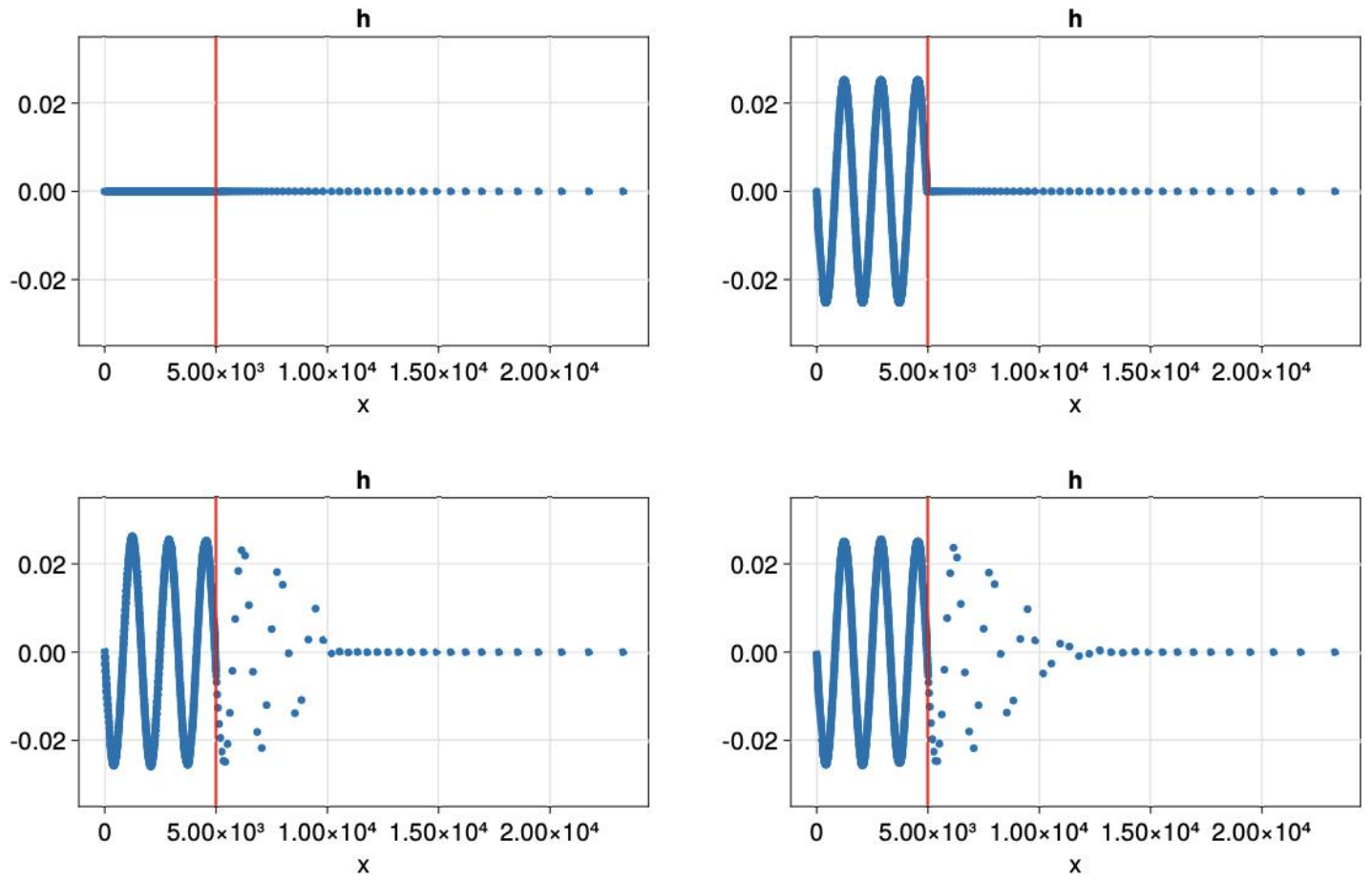}
    \caption{1D wave train: Evolution of the height perturbation $h$ over time for $x \in \Omega^F \cup \Omega^S$: Top-left: $t=0~$s, top-right: $t=500~$s, bottom-left $t=1,000~$s, bottom-right $t=5,000~$s. The waves caused by a continuous forcing at the leftmost point of the domain are transmitted through $\Omega^F$ until reaching the semi-infinite element $\Omega^S$ (the interface between the two domains is illustrated by the red line in the plots), where they are progressively damped over time. }
    \label{fig:Wave_Train}
\end{figure}

\begin{figure}[h!]
    \centering
    \includegraphics[width=0.6\textwidth]{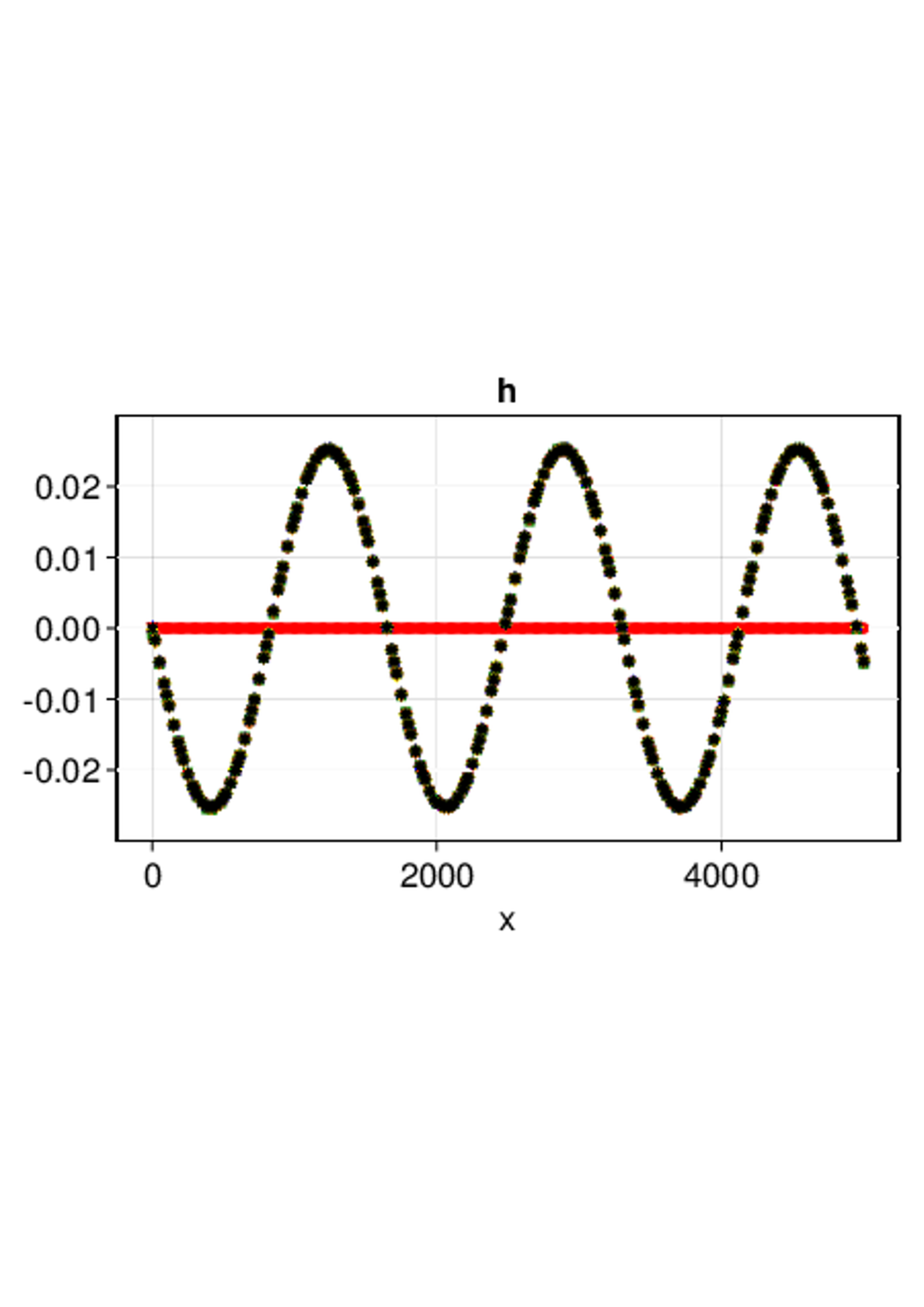}
    \caption{1D wave train: Overlapped wave train solution for $x \in \Omega^F$ \\ at times  $t \in \{0,5,10,15,20,25,30,35,40,45,50\}\times 10^3~$s. With the exception of the solution at $t=0~$s (the initial state, shown in red, is at rest), the solutions (which are periodic in time) overlap nearly perfectly demonstrating that even over long periods of time the solution is not degraded by waves reflecting off of the absorbing layer (shown in figure \ref{fig:Wave_Train}). This image was generated using a lower resolution simulation with only $N_x=50$ elements in the finite domain to allow for more easily distinguishable markers.}
    \label{fig:Wave_Train_overlap}
\end{figure}

As demonstrated in figure \ref{fig:Wave_Train}, the forcing at the boundary instigates a series of waves to move towards the end of the finite domain. Once inside the semi-infinite element, the waves are gradually damped before disappearing entirely. Observation of the finite domain shows that no noticeable changes can be seen in the wave train. To verify this, we perform a simulation of the wave train up to a final time $T_F = 50,000~$s. Figure \ref{fig:Wave_Train_overlap} shows an overlap of the solution taken at an interval $\Delta T=5,000~$s starting at $t=0~$s; This figure was generated using a low resolution simulation with only $N_x=50$ elements in the finite domain to allow the reader to better distinguish the individual points which would be, otherwise, completely overlapping and indistinguishable. The boundary condition at these times is identical and the solution should also be identical as it is periodic in time. It is clear from the figure that there is a remarkable overlap of the finite domain solutions at these times; we can thus conclude that the outgoing waves are effectively damped.

Table \ref{tab:Wave_Train} presents the time measurements for this test. Like in table \ref{tab:1D_Wave_Time}, we can see that semi-infinite elements in the absorbing layer are considerably less computationally expensive than using an extended finite domain while maintaining the same resolution and the same end point. The Laguerre approximation using elements of order 50 is nearly four times less expensive than the extended finite domain with a Rayleigh sponge. Furthermore, because the finite domain contains more elements, we can see that the contribution of the semi-infinite right hand side to the total time per time step is even less than in the previous test case for all orders of the Laguerre element that were tested. It is notable that this is the case even though the elements in the finite domain are of a lower order than the wave equation test case (4 as opposed to 6). We will show in the upcoming sections that this relative cost decreases further in 2D.

\

\begin{table}[h!]
    \centering
    \begin{tabular}{|c|c|c|c|c|c|}
    \hline
         Absorbing layer type & $T^*$ & $X_{end}$ & $T^*_{Finite} \%$& $T^*_{Laguerre} \%$ & $N_x$\\
    \hline
    \hline
         Semi-infinite element of order 15 & 1& 9984& 86.05& 13.95 & 300\\
         Semi-infinite element of order 20 & 1.05& 11837& 81.36& 18.64 & "\\
         Semi-infinite element of order 25 & 1.06& 13713& 80.55& 19.45& "\\
         Semi-infinite element of order 30 & 1.10& 15604& 78.00& 22.00& "\\
         Semi-infinite element of order 40 & 1.18& 19418& 72.70& 27.30& "\\
         Semi-infinite element of order 50 & 1.26& 23262& 68.00& 32.00& "\\
         \hline
         Extended finite domain of order 4 & 1.85& 9984& 100& N/A & 600\\
         " & 2.26& 11837& "& "& 711\\
         " & 2.65& 13713& "& "& 824\\
         " & 3.10& 15604& "& "& 937\\
         " & 3.85& 19418& "& "& 1167\\
         " & 4.63& 23262& "& "& 1398\\
        \hline
    \end{tabular}
    \caption{Timings of 1D wave train simulations with and without Laguerre semi-infinite elements in the absorbing layer. Simulations using the extended finite domain have the same physical domain size as simulations using semi-infinite elements, but the space that would be otherwise occupied by these elements is instead occupied by elements of the finite domain retaining the same resolution $\Delta x = 4.16~$m and element order $N=4$.}
    \label{tab:Wave_Train}
\end{table}

\subsection{2D advection-diffusion equation}

\begin{figure}[h!]
    \centering
	\includegraphics[width=0.99\textwidth]{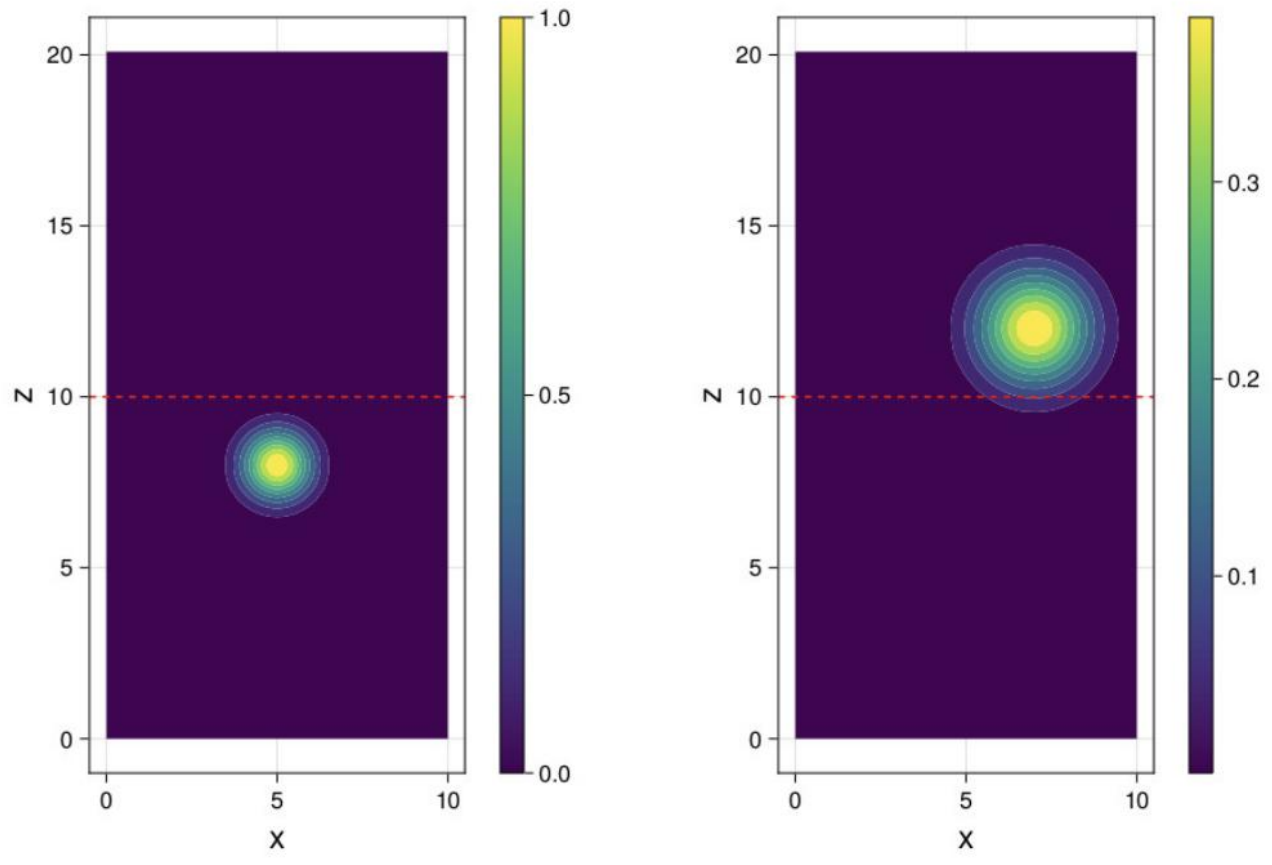}
        \caption{Solution to the 2D advection-diffusion equation. The tracer is smoothly transmitted through the $\Omega^{\rm F}$/$\Omega^{\rm S}$ interface at $z=5$ m indicated by the dashed red line.  A Laguerre semi-infinite element of order 40 is used. This test intentionally does not include an absorbing layer in order to avoid damping the solution inside $\Omega^{\rm S}$.}
      \label{fig:amr}
\end{figure}

Equations \eqref{eq:CL} with flux and source vectors \eqref{eq:2D_AdvDiff} describe the two-dimensional transport and diffusion of a passive tracer.  We performed the transport-diffusion test described in
Sec. 4.2.1 of \cite{vismaraBenacchio2023} to verify the coupling between the finite and semi-infinite domain for two dimensional cases.  The finite domain extends like $\Omega^F = [-5,5] \times [0,10]$, which is subdivided into $N_x=12$ elements in the x direction and $N_z=125$ elements in the z direction. Both directions use elements of polynomial order 4, for an effective resolution $\{\Delta x,\Delta z\} = \{0.208, 0.02\}$ m. The semi-infinite domain, defined as $\Omega^S = [-5,5]\times[10,\infty)$, is discretized using 12 Laguerre semi-infinite elements of order 40. Taking into account a scaling factor $\lambda=0.07~$m, the last LGR node on each semi-infinite element has a coordinate $Z_{end}=20.09$. To verify that the solution is not adversely affected by inter-domain interface, no damping is used for this test. We initialize the tracer field as follows:

\begin{equation}
    q(0,x,z) = \exp(-(x-x_c)^2)\exp(-(z-z_c)^2),
\end{equation}
where $x_c = 0~$m, and $z_c=8~$m. The velocity field $\mathbf{u}$ is constant and taken to be $\mathbf{u}=[0.5,1]^T {\rm m s^{-1}}$.
We use a constant viscosity coefficient $\nu = 0.1~{\rm m^2 s^{-1}}$, and a time step $\Delta t=0.0005~$s to advance the simulation to its final time $T=4~$s. Figure \ref{fig:amr} shows the solution at initial and final times and it shows that the tracer is correctly transported and diffused through the $\Omega^{\rm F}$/$\Omega^{\rm S}$, consistent with the results shown in Figure 3 in \cite{vismaraBenacchio2023}. We can see that the symmetry of the solution is maintained and no discontinuities are introduced by the interface.

\begin{figure}[h!]
    \centering
	\includegraphics[width=0.99\textwidth]{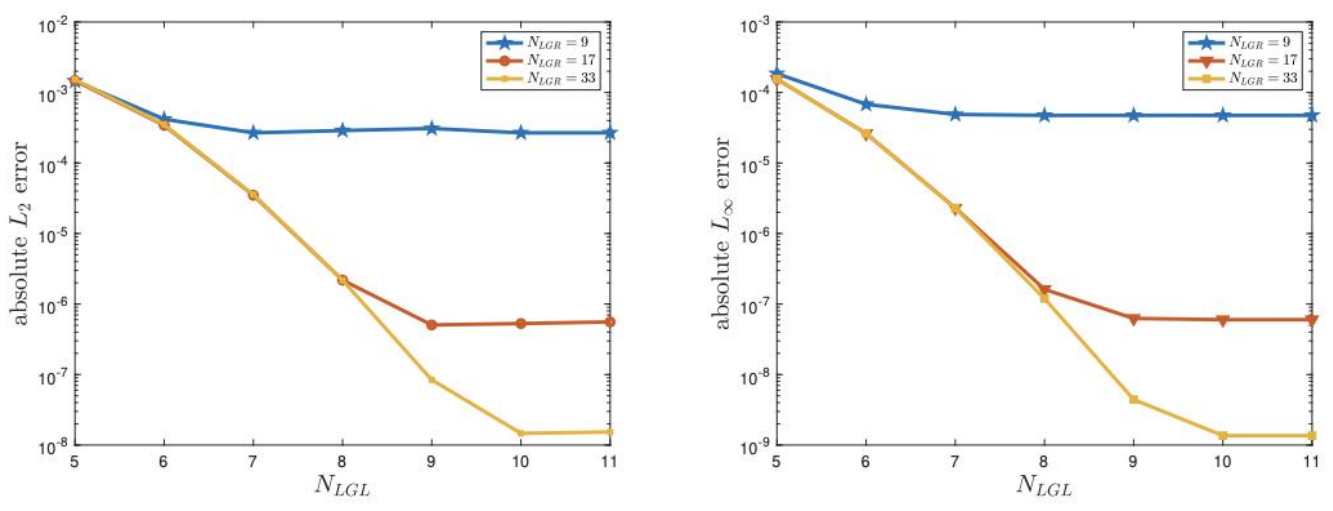}
        \caption{Absolute $L_2$ (left) and $L_{\infty}$ (right) errors as a function of $N_{LGL}$ and $N_{LGR}$ for the advection-diffusion problem.}
      \label{fig:AdvDiff_Conv}
\end{figure}

Figure \ref{fig:AdvDiff_Conv} reports the $L_2$ and $L_{\infty}$ absolute errors with respect to the exact solution as a function of the $N_{LGL}$ and $N_{LGR}$. We maintain the same number of elements as reported above and vary the order of the elements used in $\Omega^F$ and $\Omega^S$. Both panels of the figure show that the error floor is determined by the order of the semi-infinite element in the $z$ direction. Increasing the order of elements in $\Omega^F$ and in the $x$ direction of $\Omega^S$ decreases the error towards this floor. The error values we obtain are also comparable to those reported in table D1 of \cite{vismaraBenacchio2023} and this gives us confidence the accuracy of this approach.

Table \ref{tab:AdvDiff_Table} presents the timing results of this test case similarly to tables \ref{tab:1D_Wave_Time}, \ref{tab:Wave_Train}. For these tests we maintain a constant $N_x=12$, and only vary $N_z$ when no semi-infinite elements are used. This extends the finite domain while maintaining the same resolution. The table shows that the wall time of Laguerre-Legendre right hand side is considerably lower than the finite domain right hand side, costing less than a fifth of the time per time step for Laguerre semi-infinite elements of up to order 40, and less than a third for orders up to 60. Furthermore, using the Laguerre-Legendre elements is significantly cheaper than extending the finite domain with constant resolution, just as was the case for the one dimensional test cases. 

\begin{table}[h!]
    \centering
    \begin{tabular}{|c|c|c|c|c|c|}
    \hline
         Extended domain type& $T^*$ & $Z_{end}$ & $T^*_{Finite} \%$& $T^*_{Laguerre} \%$ &$N_z$\\
    \hline         
    \hline
         Semi-infinite elements of order 15 & 1& 13.4& 93.11& 6.89& 125\\
         Semi-infinite elements of order 20 & 1.02& 14.7& 91.02& 8.98& "\\
         Semi-infinite elements of order 30 & 1.06& 17.4& 87.20& 12.80& "\\
         Semi-infinite elements of order 40 & 1.12& 20.09& 82.96& 17.04& "\\
         Semi-infinite elements of order 50 & 1.17& 22.7& 79.20& 21.80& "\\
         Semi-infinite elements of order 60 & 1.24& 25.4& 74.51& 25.49& "\\
         \hline
         Extended finite domain of order 4 & 1.25& 13.4& 100& N/A& 167\\
         " & 1.38& 14.7& "& "& 183\\
         " & 1.64& 17.4& "& "& 217\\
         " & 1.92& 20.09& "& "& 251\\
         " & 2.21& 22.7& "& "& 284\\
         " & 2.54& 25.4& "& "& 317\\
        \hline
    \end{tabular}
    \caption{Timings of 2D advection-diffusion simulations with and without Laguerre semi-infinite elements. Simulations using the extended finite domain have the same physical domain size as simulations using semi-infinite elements, but the space that would be otherwise occupied by these elements is instead occupied by elements of the finite domain. In this case the number of elements in the $x$-direction $N_x$ remains the same but the number of elements in the $z$-direction $N_z$ is adjusted. This is done so that enough elements of order $N=4$ are used to maintain a constant vertical resolution $\Delta z= 0.02~$m.}
    \label{tab:AdvDiff_Table}
\end{table}

\subsection{2D Helmholtz equation}
We Helmholtz equation is obtained from (\ref{eq:CL}) by means of vectors (\ref{eq:Helmholtz}). We solved it on the 2D semi-infinite channel $\Omega = [0, \infty) \times [-\pi/2, \pi/2]$ subject to a zero Dirichlet boundary conditions (BCs) on the boundary.  This test verifies the convergence rate of the solution and presents some guidance on selecting the orders of the elements in the finite and semi-infinite domains. Applying nodal SEM with inexact integration to this equation along with the Dirichlet BCs yields the matrix-vector problem
\begin{equation}
-\mathbf{L}_{IJ} u_J + \alpha^2 \mathbf{M}_{IJ} u_J = -\mathbf{M}_{IJ} f_J
\label{matvec}
\end{equation}
where $\mathbf{L}_{IJ}$ is the block diagonal Laplacian matrix (see Sec. 8.6.4 in \cite{giraldoBOOK}) and $\mathbf{M}_{IJ}$ is the diagonal mass matrix given by \eqref{eq:matvec3}.  The matrix-vector problem \eqref{matvec} is solved via a direct sparse solve by considering the exact solution
\begin{equation}
u(x,y) = e^{-x/L} \sin\left( \frac{x}{L} \right)\cos(y),
\label{sol}
\end{equation}
with $L=2$ m and then analytically construct the RHS $f(x,y)$ using the method of manufactured solutions using a wavenumber $\alpha  = 10 \mbox{m}^{-1}$.  

\begin{figure}[h!]
\centering
\includegraphics[width=0.9\textwidth]{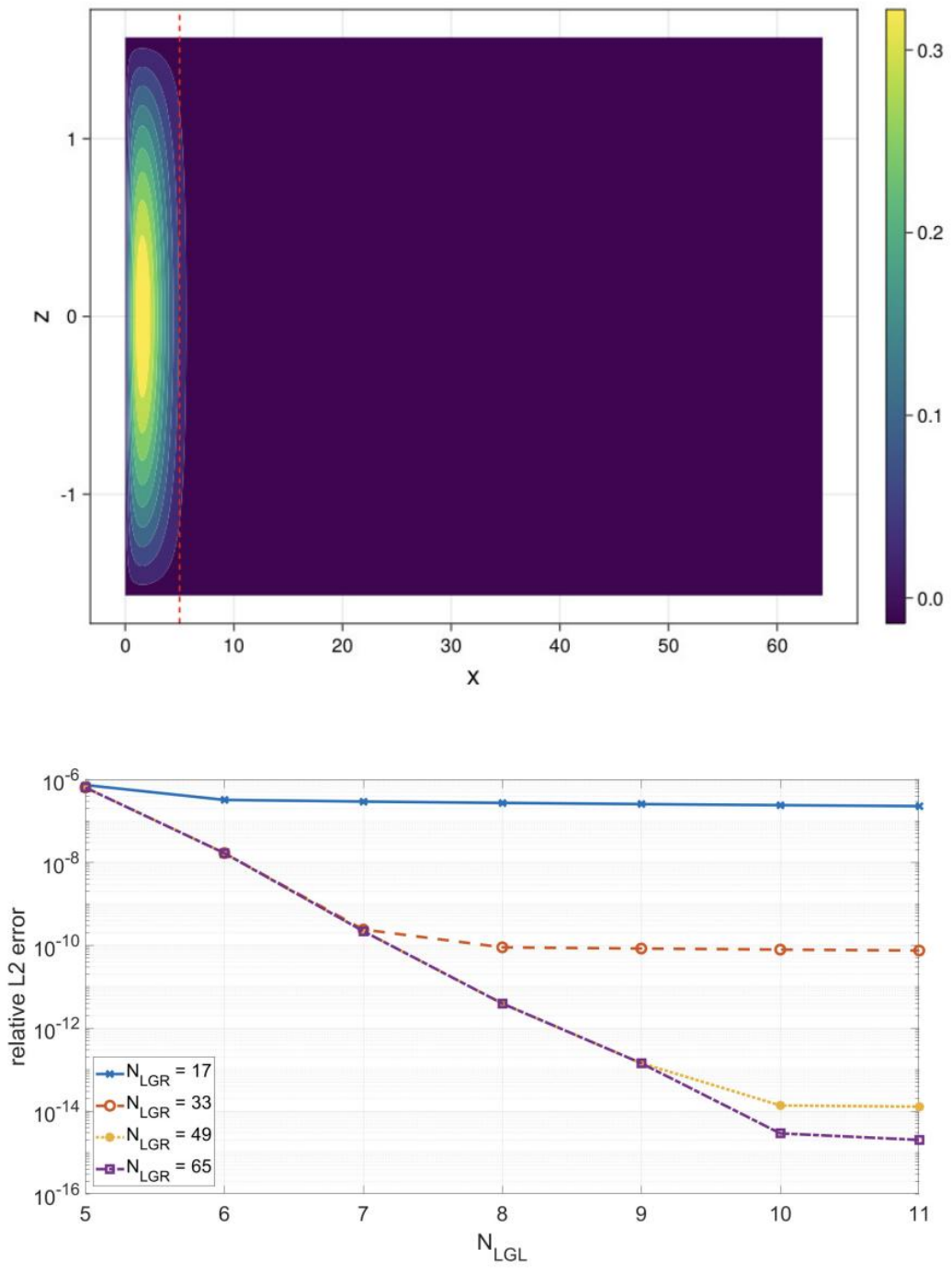}
\caption{ Top: Solution to the Helmholtz equation with $N_{LGL} = 11$ and $N_{LGR} = 65$, the red line shows the location of the interface between $\Omega^F$ and $\Omega^S$. Bottom: relative $L^2$ 
error as a function of $N_{LGL}$ and $N_{LGR}$.}
\label{fig:helmholtz} 
\end{figure}

We partitioned the channel into a finite domain $\Omega^F = [0,5] \times [-\pi/2,\pi/2]$ and semi-infinite domain $\Omega^S = [5,\infty) \times [-\pi/2,\pi/2]$.  The domain $\Omega^F$ was discretized with $N_F = 16$ finite spectral elements with orders ranging from 4 to 10, while $\Omega^S$ was discretized with $N_S = 4$ semi-infinite elements with orders ranging from 16 to 64.  We then computed the relative $L^2$ error of the numerical solution relative to the analytical solution \eqref{sol}.  Figure \ref{fig:helmholtz} displays a) the numerical solution to the Helmholtz equation with $N_{LGL} = 11$ and $N_{LGR} = 65$, b) relative $L^2$ error as a function of $N_{LGL}$ and $N_{LGR}$.   

Close examination of panel b) shows that the choice of the semi-infinite element order $N_{LGR}$ determines a floor for the relative error, while the increase of $N_{LGL}$ reduces the relative error to this floor.  In particular, choosing $N_{LGR} = 48$ and $N_{LGL} = 9$ yields
a relative error of $\sim 10^{-14}$.  Since this problem has a wavelike solution and an exact solution, it provides guidance on choosing the order of the LGL and LGR nodal basis functions. The next test case will show that semi-infinite elements can be used in the absorbing layer for reducing the cost in a benchmark atmospheric flow test case.

\subsection{2D Euler equations}
\paragraph{Rising thermal bubble}
We validate the Laguerre-Legendre semi-infinite element approach on the compressible Euler equations by using it to simulate a classic rising thermal bubble case. In this test a perturbation $\Delta \theta$ is introduced to a neutral atmosphere with uniform potential temperature $\theta_0 = 300~$K. The finite domain is $\Omega^F = [-5,5]~{\rm km} \times [0,5]~{\rm km}$, which is subdivided into $N_z=N_x=20$ elements of order $N=4$. This yields a horizontal resolution $\Delta x = 125~$m and a vertical resolution $\Delta z = 62.5~$m. The potential temperature perturbation and initial pressure are defined as follows:

\begin{subequations}
\label{eq:rtb}
\begin{equation}
    \Delta \theta = \theta_c \left(1-\frac{r}{r_0}\right) \mbox{  if } |r| \leq r_0
\end{equation}
\begin{equation}
    r = \sqrt{(x-x_c)^2 + (z-z_c)^2} 
  \end{equation}
\begin{equation}
    p = p_0 \left(1-\frac{gz}{c_p \theta} \right)^{cp/R},
  \end{equation}
\end{subequations}
where $\theta_c = 2~$K, $r_0 = 2000~$m, $x_c=0~$m, $z_c=2500~$m, and $p_0=1000~$hPa. The initial potential temperature is $\theta = \theta_0 + \Delta\theta$, and the initial density can be deduced through the ideal gas law \eqref{eq:eos}. 

A set of 20 Laguerre-Legendre semi-infinite elements of order 24 is added on top of $\Omega^F$ to build the semi-infinite domain $\Omega^S=[-5,5]~{\rm km} \times [5,\infty)~{\rm km}$. To verify that this approach is able to solve the Euler equations without adversely affecting the solution, we do not use a damping layer for this test and simply verify that the rising thermal crosses the $\Omega^F-\Omega^S$ boundary without being affected by the change in element type and drastic change in resolution. To stabilize the solution past $t=500~s$ a viscosity coefficient $\nu = 30~{\rm m^2s^{-1}}$ is used in conjunction with a constant thermal diffusivity $\kappa = 2\nu$. This is standard and would be required even for simulations without a semi-infinite element. Free-slip type boundary conditions are used at all domain edges.

\begin{figure}
    \centering
    \includegraphics[width=0.85\textwidth]{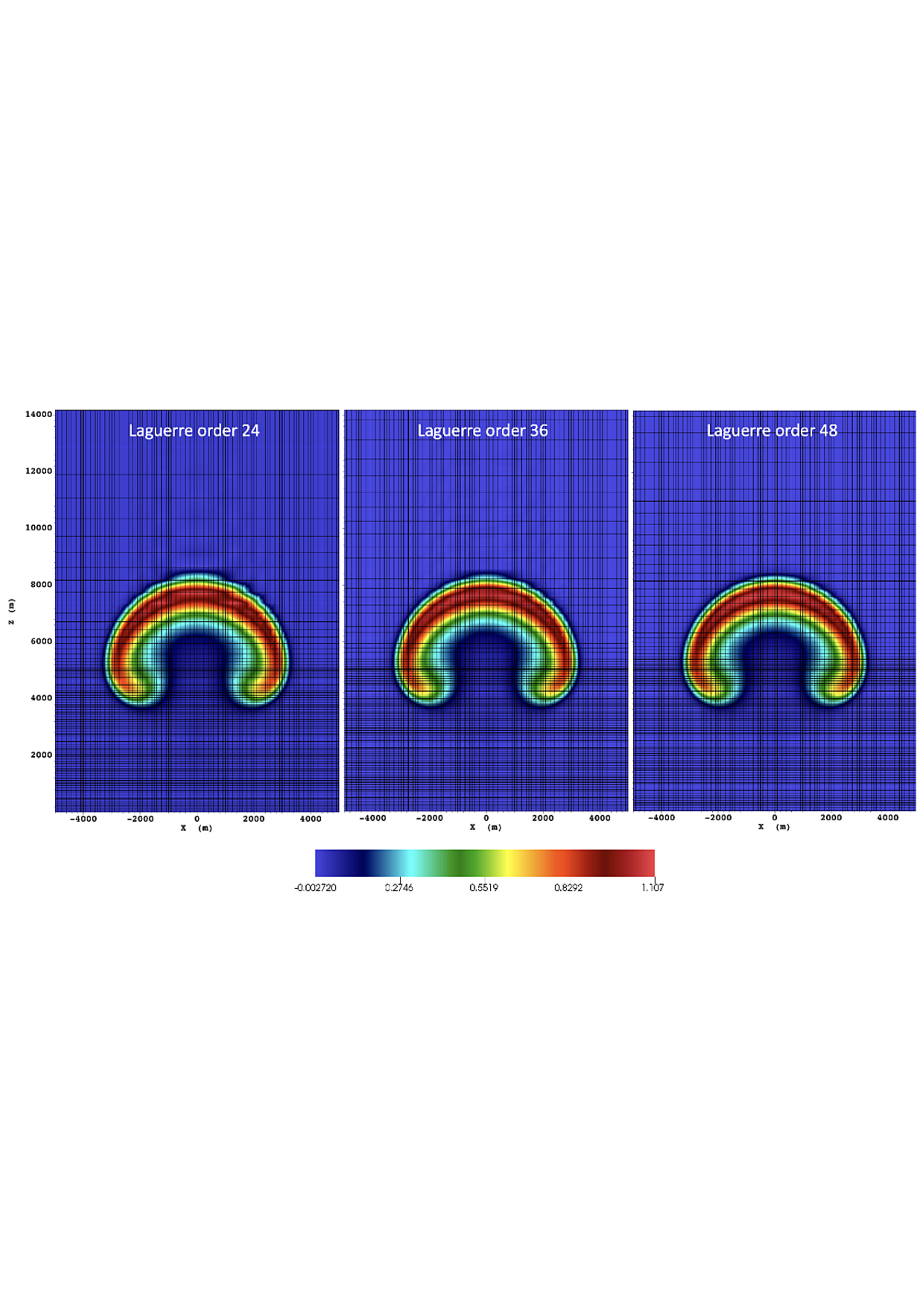}
    \caption{Potential temperature perturbations for the 
    rising thermal bubble at $t=1000~$s using semi-infinite elements of order 24 (left), 36 (middle), and 48 (right). In all three cases the bubble rises correctly and is able to pass through the interface between $\Omega^F$ and $\Omega^S$ at $z=5$ km without any issues. Note that a higher order semi-infinite element yields a sharper, higher resolution solution. The grid used for this test is fully visible to visually show the transition between the two domains. Note that no absorbing layer is used in this test since the goal is not that of testing non-reflecting conditions, but to assess the sanity of the Legendre-Lagerre discretization.}
    \label{fig:bubble}
\end{figure}

The initial potential temperature perturbation generates positive buoyancy, which causes the bubble to rise. Figure \ref{fig:bubble} shows the solution at $t=1000~$s after a significant portion of it has transported through the interface between $\Omega^F$ and $\Omega^S$ at $z=5$ km. The figure shows that the solution is correctly transported and diffused through the interface. The symmetry of the bubble in the $x$ direction is maintained, and the interface does not introduce spurious noise or discontinuities. The figure also highlights the improvement in solution quality as the order of the semi-infinite element is increased. The sharpest solution is obtained using a semi-infinite element of order 48 which is displayed in the right panel of figure \ref{fig:bubble}. With this approach validated for the Euler equations, we move on to demonstrating that we can use it for effective absorbing layers in atmospheric flows.

\paragraph{Linear hydrostatic mountain waves}
In this test, a constant horizontal flow with velocity $U=20~\rm{ms^{-1}}$ impinges on a mountain in a stratified atmosphere. These flow conditions and the mountain's size determines the structure of the resulting waves.
As soon as the flow encounters the mountain, gravity waves quickly propagate both horizontally and vertically.  In a proper implementation, these waves should leave the domain with no reflection.

The background state for this test is a hydrostatically balanced atmosphere whose pressure and potential temperature are:
\begin{subequations}
\label{eq:rtbini}
\begin{equation}
p = p_0 \left[ 1 + \frac{g^2}{c_p \theta_0 N^2} \left(\exp\left(\frac{-z \mathcal{N}^2}{g} \right) - 1\right)\right]^{c_p/R}
\end{equation}
\begin{equation}
\theta = \theta_0 \exp\left(\frac{z \mathcal{N}^2}{g}\right),
\end{equation}
\end{subequations}
where $p_0=1000~$hPa and $\theta_0=250~$K are the sea level values of pressure and potential temperature, and 
\[
\mathcal{N} = \frac{g}{\sqrt{c_p \theta_0}}=0.0196~{\rm s^{-1}} 
\] is the Brunt-V\"ais\"al\"a frequency.
An Agnesi mountain with height $h=1~$m and half-width $a=10,000~$m is located at the center of the domain $x_c=0~$m with shape
\begin{equation}\label{agnesi mountain}
    z = \frac{h a^2}{(x - x_c)^2 + a^2}.
\end{equation}

For this test as well as for the other case involving topography, a terrain following sigma coordinate \cite{galChenSomerville1975a} is used.
We consider a finite domain $\Omega^F = [-120, 120]~{\rm km} \times [0,15]~$km, and subdivided it into $N_x \times N_z= 120 \times 21$ elements of order 4, which lead to the resolution $(\Delta x, \Delta z) = (500~{\rm m},178~{\rm m})$. A set of 120 Laguerre-Legendre semi-infinite elements are added on top of $\Omega^F$ yielding the semi-infinite domain $\Omega^S=[-120, 120]~{\rm km} \times [15,\infty)~$km. In the horizontal direction, each semi-infinite element uses an order 4 spectral element discretization using LGL nodes, and in the vertical direction it uses an order 14 Laguerre function basis on LGR nodes and a scaling factor $\lambda=300~$m which yields an end point $Z_{end} = 28853$. Due to prevalence of both acoustic waves and gravity waves, as well as the sensitivity of high order numerical methods, we make use of the spectral low pass Boyd-Vandeven filter \cite{boyd1996erfc} to help insure the stability of the simulation (this approach is standard for spectral element simulations of the atmosphere). The filter on the elements in $\Omega^S$ is different from that in $\Omega^F$ and the continuity of the solution at shared interface nodes is insured by the DSS operation.

We use a sine squared function to define the damping coefficient in the vertical absorbing layer:

\begin{equation}
    \gamma (z) = \Delta \gamma \sin^2\left(\frac{\pi}{2}\frac{z-z_s}{z_{max}-z_s}\right),
\end{equation}
where $z_s=15 000$ m, $\Delta \gamma = 0.1 \mbox{s}^{-1}$ and $z_{max}$ is the top of the absorbing layer. The lateral boundaries are periodic and lateral damping layers are used on each side (these damping layers are in $\Omega^F$ and are not on semi-infinite elements). These boundary conditions are also used for the remaining other mountain test presented in this paper. 

\begin{figure}[h!]
\centering
\includegraphics[width=0.8\textwidth]{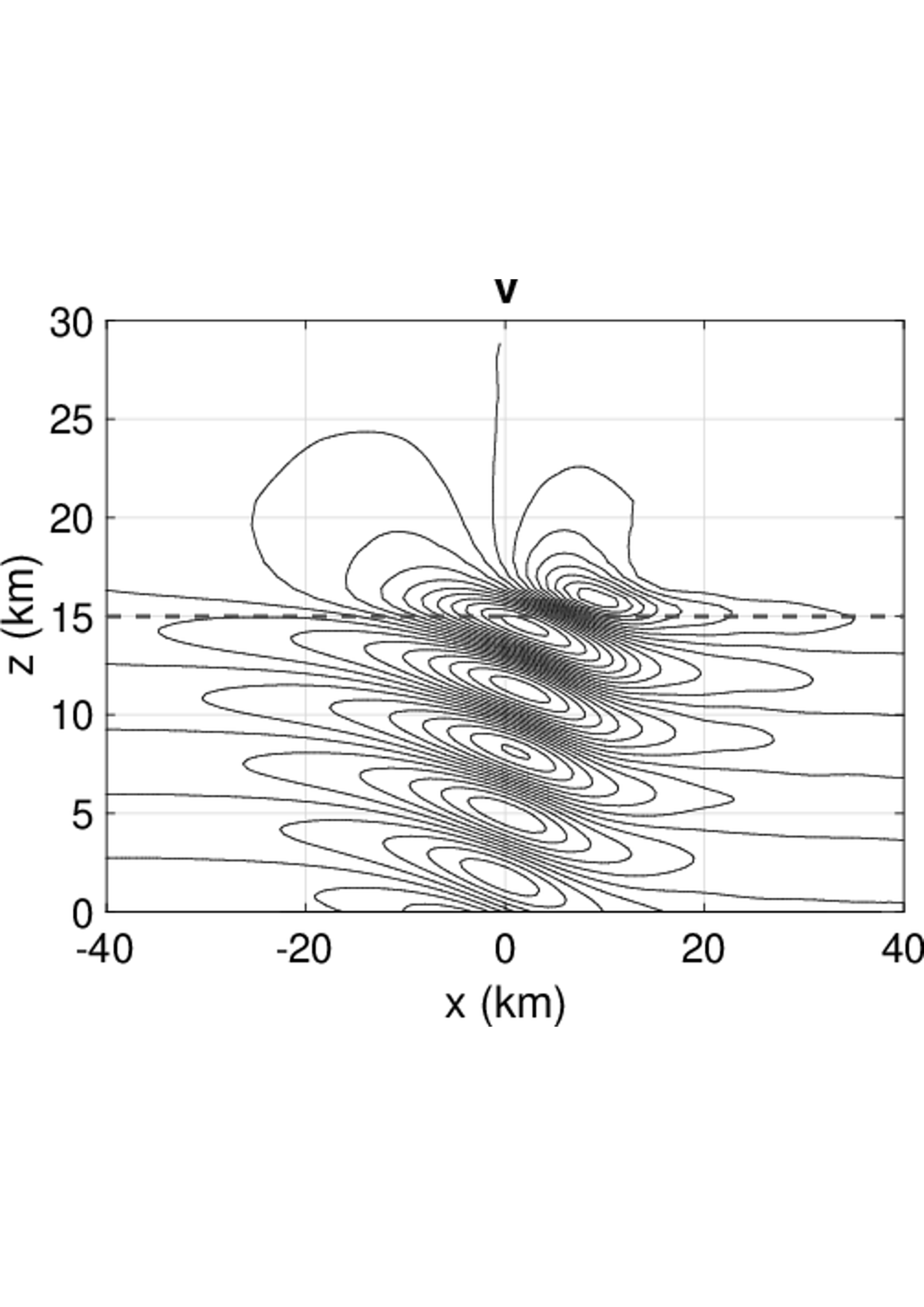}
\caption{Time-converged numerical vertical velocity $w$ for the
Linear Hydrostatic Mountain.  The numerical solutions are displayed at t=30,000 s for $\Delta x = 500~$m, $\Delta z = 178~$m. The contours of the vertical velocity range between -0.005 and +0.005 with interval 0.0005. The dashed line shows where $\Omega^F$ meets $\Omega^S$.
In $\Omega^S$ there is only one element of order 14 along the vertical direction (see Fig.~\ref{fig:twodomains} for reference.)}
\label{fig:mountain_uw}
\end{figure}

Figure \ref{fig:mountain_uw} shows the contours of the vertical velocity at $t=30,000~$s. The solution is stable, the outgoing waves are effectively damped within $\Omega^S$, and the solution is physically meaningful when compared to other numerical solutions of atmospheric models using spectral elements see e.g. \cite{sridharEtAl2019CLIMALES}.

Table \ref{tab:hydrostatic_table} shows the time per time step for different configurations of the linear hydrostatic mountain test case. A fixed domain end point $Z_{end} = 30,000~$m  is maintained for the cases where no semi-infinite element is used and instead only the number of vertical elements is changed. Even while only using ten additional order 4 spectral elements in the vertical direction to replace the Laguerre semi-infinite element the time per time step remains about 10\% higher than running a single order 20 Laguerre element. If we would seek to run the entire vertical domain at the same resolution as the finite domain of these simulations utilizing a Laguerre semi-infinite element, the cost increase would be nearly 50\% higher. Considerably few additional vertical spectral elements (we estimate no more than three) would have to be used to obtain a similar time per time step as using a single semi-infinite element in this case and this would also yield a much lower order discretization within the absorbing layer. 

\begin{table}[h!]
    \centering
    \begin{tabular}{|c|c|c|c|c|c|}
    \hline
         Absorbing layer type& $T^*$ & $Z_{end}$ & $T^*_{Finite}$ \%& $T^*_{Laguerre}$ \%&$N_z$\\
    \hline
    \hline
        Semi-infinite elements of order 14 & 1 & 28853 & 90.41 & 9.59 & 21\\
        Semi-infinite elements of order 18 & 1.02 & 33279 & 87.97 & 12.03 & 21\\
        Semi-infinite elements of order 20 & 1.04 & 35513 & 86.49 & 13.51 & 21\\
        \hline
        Extended finite domain of order 4 & 1.14 & 30000 & 100 & N/A & 31\\
        " & 1.29 & " & " & " & 35\\
        " & 1.47 & " & " & " & 40\\
        \hline
    \end{tabular}
     \caption{Timings of linear hydrostatic mountain simulations with and without Laguerre semi-infinite elements in the absorbing layer. A fixed vertical domain end point is used for simulation using only a finite domain while the simulations using semi-infinite elements are allowed to have varying $Z_{end}$ depending on the order of the semi-infinite elements. In this case the number of elements in the $x$-direction $N_x$ remains the same but the number of elements in the $z$-direction $N_z$ is adjusted. The additional elements are all of order $N=4$ but the vertical resolution $\Delta z$ is allowed to change.}
     \label{tab:hydrostatic_table}
\end{table}

\begin{figure}
   \centering
    \includegraphics[width=0.8\textwidth]{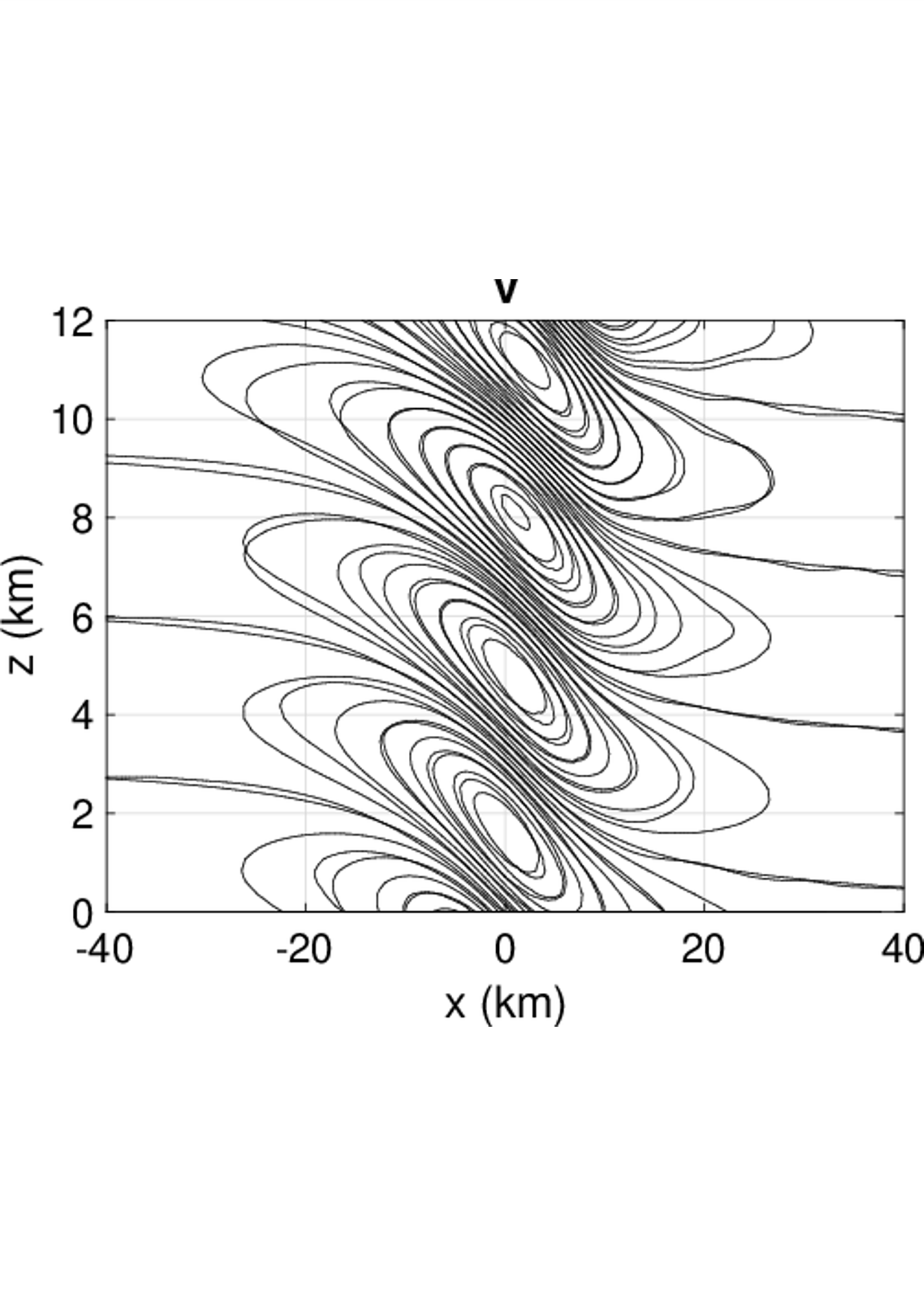}
    \caption{Comparison of vertical velocity contours of the numerical solution (blue) of the linear hydrostatic mountain wave problem with an approximate linear Fourier solution (red). The contours are the same contours used in the vertical velocity plot of figure \ref{fig:mountain_uw}.  The numerical solutions compare favorably to the analytical approximation and display the expected vertical wavelength $\lambda_z = 2 \pi U / \mathcal{N} \sim 6.4$ km. The Fourier solution is found under specific assumptions that are not met by the numerical model. For this reason, the small differences between the two contour sets are expected.}
     
    \label{fig:Mountain_overlap}
\end{figure}

In order to complete the validation of the hydrostatic and mountain wave simulations, we constructed a linear Fourier solution using the approach outlined in \cite{smith1980linear}. The vertical velocity is expressed as a Fourier integral under
the Boussinesq approximation that is then evaluated using adaptive Gauss quadrature in wavenumber space.  Since this solution neglects vertical variations in density, it is only valid for heights less than a scale height ($\sim$9 km). Figure \ref{fig:Mountain_overlap} compares vertical velocity profiles of the numerically obtained solution (in blue) with the linear Fourier solution (in red). The figure shows a good overlap of the two solutions and the existing deviations are expected. This is given that the two models deviate due to the analytical solution relying on the linearization of the Euler equations and the Boussinesq approximation. Furthermore, this deviation is comparable to other atmospheric models.

\paragraph{Sch\"ar mountain waves}
The Sch\"ar test \cite{scharLeuenberger2002} consists of a uniform flow with a reference horizontal velocity $U=10~\rm{ms^{-1}}$  in a stratified atmosphere with Brunt-V\"ais\"al\"a frequency $\mathcal{N}=0.01$~$\rm{s^{-1}}$, sea-level pressure $p_0=1000~$hPa and and potential temperature at sea-level $\theta_0=280~$K. 
The flow impinges the five-peak mountain defined as 
\[
z = h e^{-\left(\frac{x}{a} \right)^2}\cos^2\left(\frac{\pi x}{\lambda_c} \right)
\]
with parameters $h=250$~m, $a=5000$~m, and $\lambda_c=4000$~m.

We consider a finite domain $\Omega^F = [-25, 25]~{\rm km} \times [0,15]~$km, and discretize it using $N_x=20$ elements in the x direction and $N_z=7$ elements in the z direction. These elements are of polynomial order 10 in both directions, leading to an effective resolution $(\Delta x, \Delta z) = (214~{\rm m},250~{\rm m})$. A set of 20 Laguerre-Legendre semi-infinite elements are added on top of $\Omega^F$ which translates to the semi-infinite domain $\Omega^S=[-25, 25]~{\rm km} \times [15,\infty)~$km. In the horizontal direction, each semi-infinite element uses an order 10 spectral element discretization using LGL nodes, and in the vertical direction it uses an order 14 Laguerre basis function on LGR nodes and a scaling factor $\lambda=300~$m which yields an end point $Z_{end} = 28853$ m. Similarly to the previous mountain wave tests a Boyd-Vandeven filter \cite{boyd1996erfc} is used to help insure the stability of the simulations and the a sine squared function is used in the damping coefficient of the Rayleigh damping layer.
\begin{figure}
\centering
\includegraphics[width=0.8\textwidth]{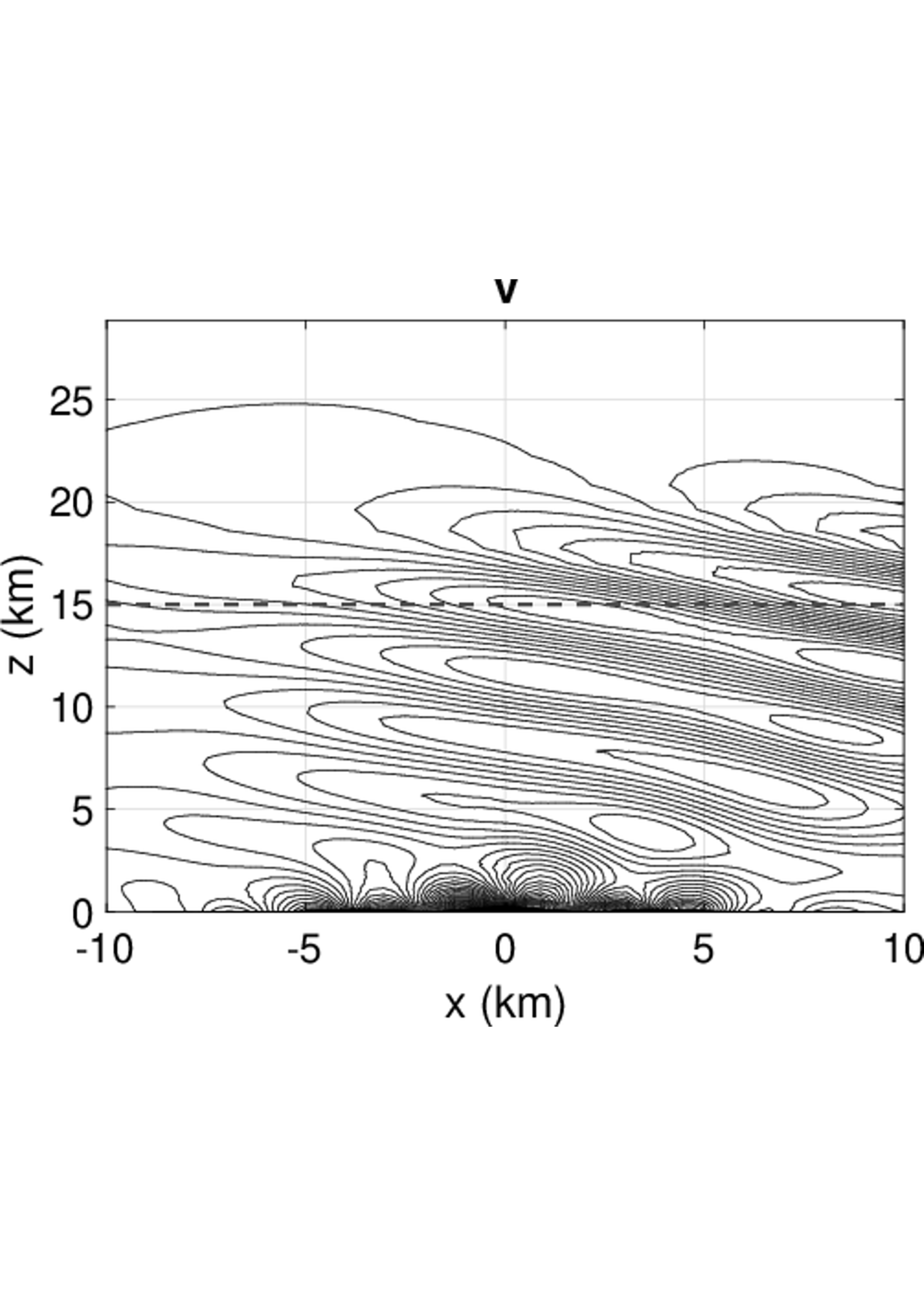}
\caption{Sch\"ar mountain waves. Time-converged vertical velocity at t=36,000 s for $\Delta x = 250~$m, $\Delta z = 220~$m.
The contours of the vertical velocity range between $-2.0~\rm{ms^{-1}}$ and $+2.0~\rm{ms^{-1}}$ with interval $0.1~\rm{ms^{-1}}.$ The dashed line shows where $\Omega^F$ meets $\Omega^S$.
In $\Omega^S$ there is only one element of order 14 along the vertical direction (see Fig.~\ref{fig:twodomains} for reference.)}
\label{fig:schar_mount}
\end{figure}

\begin{figure}
\centering
\includegraphics[width=0.8\textwidth]{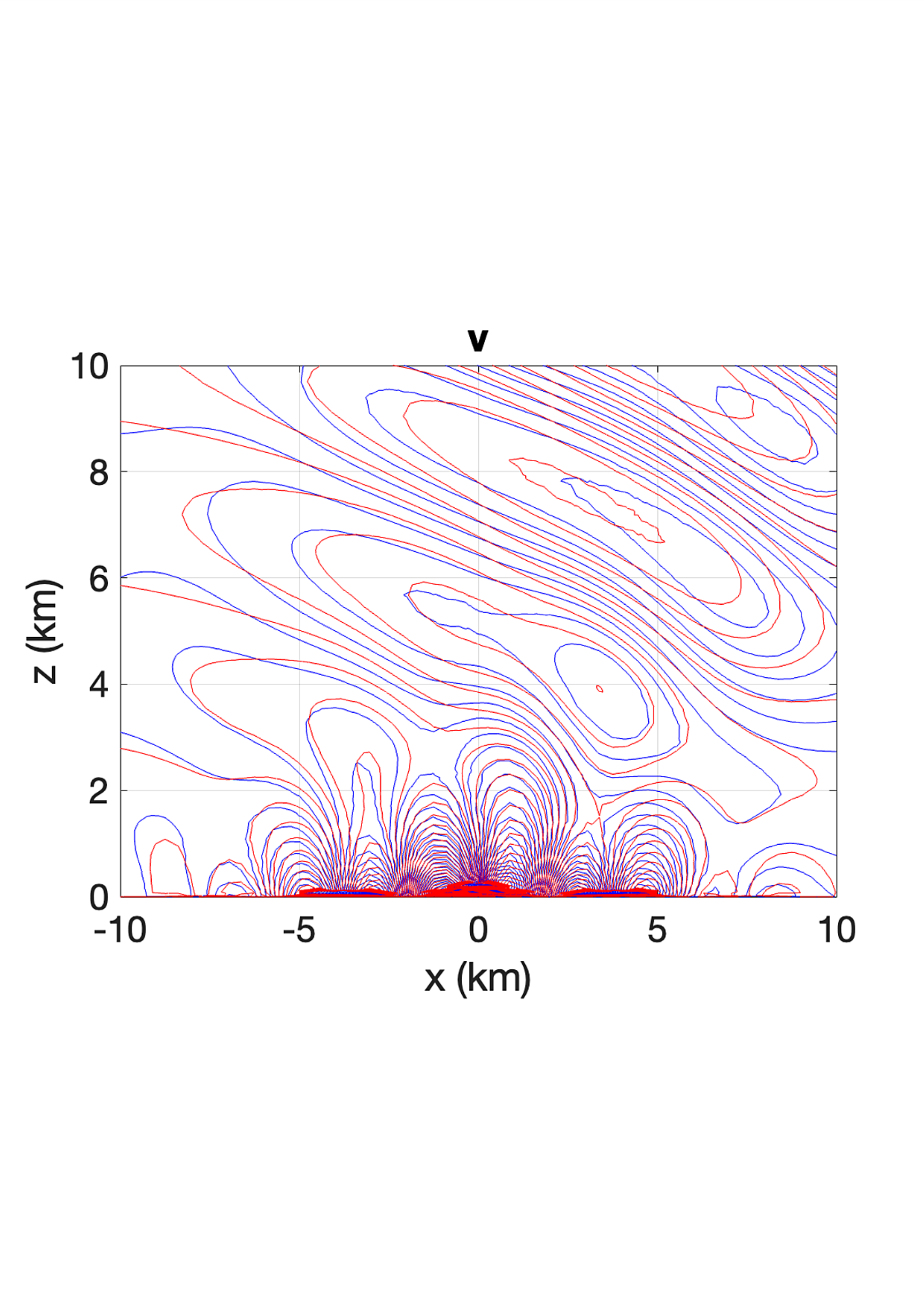}
\caption{Comparison of vertical velocity contours of the numerical solution (blue) and semi-analytical Fourier solution (red) for the Sch\"ar mountain test. As for those shown in Fig.~\ref{fig:Mountain_overlap}, some discrepancy is due to the linearization of the Euler equations used to calculate the semi-analytic Fourier solution.}
\label{fig:schar_overlap}
\end{figure}

\begin{table}[h!]
    \centering
    \begin{tabular}{|c|c|c|c|c|c|}
    \hline
         Absorbing layer type& $T^*$ & $Z_{end}$ & $T^*_{Finite}$ \%& $T^*_{Laguerre}$ \%&$N_z$\\
    \hline
    \hline
        Semi-infinite elements of order 14 & 1 & 28853 & 88.62 & 11.38 & 7\\
        Semi-infinite elements of order 18 & 1.03 & 33279 & 86.02 & 13.98 & 7\\
        \hline
        Extended finite domain of order 10 & 0.92 & 21000 & 100 & N/A & 8\\
        " & 1.03 & " & " & " & 9\\
        " & 1.14 & " & " & " & 10\\
        " & 1.27 & " & " & " & 11\\
        \hline
    \end{tabular}
     \caption{Timings of Sch\"ar mountain wave simulations with and without Laguerre semi-infinite elements in the absorbing layer. A fixed vertical domain end point is used for simulation using only a finite domain while the simulations using semi-infinite elements are allowed to have varying $Z_{end}$ depending on the order of the semi-infinite elements. In this case the number of elements in the $x$-direction $N_x$ remains the same but the number of elements in the $z$-direction $N_z$ is adjusted. The additional elements are all of order $N=10$ but the vertical resolution $\Delta z$ is allowed to change.}
     \label{tab:schar_table}
\end{table}

Figure \ref{fig:schar_mount} shows the vertical velocity of the time converged numerical solution at $t=36,000~$s. The solution is stable and outgoing waves are effectively damped withing $\Omega^S$. The solution is also physically meaningful and comparable to other numerical solutions of atmospheric models, see e.g \cite{scharLeuenberger2002}. Figure \ref{fig:schar_overlap} shows an overlap of the numerical vertical velocity with a linear Fourier solution under the anelastic approximation \cite[Eq.(A10)]{klemp2015}. The figure shows a good overlap of the two solution and the existing deviations are within expectations given the differences between the models. Furthermore, these deviations are comparable to other atmospheric models.

Table \ref{tab:schar_table} shows the time per time step for different configurations of the sch\"ar mountain problem. A fixed domain end point $Z_{end}=21,000~$m is maintained for the cases where no semi-infinite element is used and instead only the number of vertical elements is changed. For this case, even while using only 4 additional elements to replace the Laguerre semi-infinite element and resulting in nearly equivalent resolutions within the finite domain, the cost per time step is nearly $25\%$ higher than using semi-infinite elements of order 18 while also having a smaller domain height. Only two additional elements of order 10 can be used to replace the semi-infinite elements if a similar time to solution is desired, but this will come at the cost of the order of accuracy within the damping layer. 

\emph{\textbf{Remark:}} For the two mountain wave test cases, we have compared solutions obtained using a standard spectral element approach against solutions obtained using the proposed approach of using semi-infinite elements in the absorbing layer. We report the root mean squared error resulting from these comparisons in Appendix C and conclude that the proposed semi-infinite element approach does not introduce any significant errors relative to a standard spectral element approach. 
\section{Discussion}
\label{sec:disc}

\subsection{Mass conservation}
As discussed in \cite{thuburn2008some}, mass conservation is a key ingredient to ensure accurate forecasts, such as correctly diagnosing surface pressure.  A gradual, secular loss of mass may distort the resolved dynamics over long term time-integration.  In addition, accurate surface pressure prediction depends on a proper mass budget \cite{thuburn2008some}.  Also, mass conservation is a prerequisite for conserving other invariants, such as momentum and total energy.  Many spectral element atmospheric models, including HOMME \cite{taylorfournier2010} and NUMA \cite{kellyGiraldo2012}, have been designed explicitly to conserve mass. 
Table \ref{tab:mass_loss} presents the relative mass loss for the method with respect to the compressible Euler equations and its sensitivity to different time integration schemes. These tests were done on a rising-thermal bubble and demonstrate that the method can indeed conserve mass.

\begin{table}[h!]
\centering
{\begin{tabular}{lcccc|cc}
 Time integrator && &  &  & & Mass loss\\ \midrule
 SSPRK54 &  & & &  &  & 1.9545155453050947e-16\\ 
 SSPRK33 &  & & &  & & 1.1727093271830568e-15\\ 
 MSRK5  &  & & &   & & 3.9090310906101895e-16\\
 \bottomrule
\end{tabular}}
\caption{Relative mass loss for the compressible Euler equations and sensitivity to the time integrators of \textit{DifferentialEquations.jl} for a rising-thermal bubble test. Mass conservation measured at t = 1000 s with active viscosity. Results are for \textbf{inexact} integration.}
\label{tab:mass_loss}
\end{table}

The compressible Euler equations, as formulated by \eqref{eq:2D_Euler},
are mass conserving since a) the continuity equation is given in conservation form, b) the continuity equation is discretized in a weak form, and c) a no-mass flux, or rigid BC, is utilized at the lower boundary.  However, since an absorbing layer is applied to the continuity equation, the model will not conserve mass. This is a common failing of explicit Rayleigh damping approaches.  To remedy this problem, implicit absorbing layers  \cite{klempDudhia2008} that only damp  vertical momentum without affecting other prognostic variables, may be applied.  
In the future, we may consider using this type of implicit sponge to allow for mass conservation.

\subsection{High-Altitude simulations}
High-Altitude (HA) weather simulations extend terrestrial weather prediction into the thermosphere, which exists from approximately 100 km to 600 km above ground level.  The thermosphere is a hot, rarefied layer of the atmosphere which contains a charged layer known as the ionosphere.  The U.S. Navy is currently developing a HA thermosphere-ionospheric model based on the pre-operational NEPTUNE model \cite{reinecke2016}.  The proposed semi-infinite element approach is a natural solution for preventing reflection of vertically propagating gravity waves.  However, there are several key differences between HA weather and terrestrial weather that need to be addressed. 

As a prototype problem, consider a 1D transport-diffusion equation for temperature $T(t,z)$:
\begin{equation}
\label{eq:1dtemp}
\frac{\partial T}{\partial t} + v \frac{\partial T}{\partial z} = \frac{\partial}{\partial z} \left( \nu \frac{\partial T}{\partial z} \right)
\end{equation}
where $v$ is vertical velocity and $\nu$ is the coefficient of kinematic viscosity.
Since $\nu$ is inversely proportional to density, the diffusive tendency grows exponentially with respect to height $z$ at high altitudes.  This physical diffusion provides a ``natural" sponge that absorbs vertically propagating gravity (and acoustic) waves.  However, the temperature at the top of the thermosphere, or the exospheric temperature $T_{\infty}$, is non-zero and undergoes large day to night swings \cite{hargreaves1992}.  Hence, the prototype problem \eqref{eq:1dtemp} does not respect the exponential decay requirement of the semi-infinite element approach.

To address this problem, we may define an auxiliary variable $\tilde{T} (t,z) = T(t,z) - T_{\infty} (t)$ that incorporates the exospheric temperature.  The auxiliary temperature $\tilde{T} \rightarrow 0$ as $z \rightarrow \infty$, allowing us to apply the semi-infinite element approach without resorting to an absorbing layer.  However, the lower boundary condition now becomes time-dependent.  We may use this approach to replacing the existing rigid upper BC/sponge layer in NEPTUNE.

\subsection{Dimensions and complexity}
Since efficiency and simulation speed is a significant concern when performing numerical simulations, we analyze the complexity of the method for up to three dimensions.

The computational cost of spectral element methods may be reduced by 1) using a tensor-product basis and 2) evaluating derivative operations using \emph{sum factorization} \cite[Chapter 4]{deville2002high}. Using these two strategies, the cost of evaluating a derivative in a $d$-dimensional domain along direction $1 \leq i \leq d$ is $\mathcal{O} \left( N_e p_1 p_2 \cdots p_i^2 \cdots p_d \right)$.  In the standard SEM, where the same basis functions are used along each direction, this reduces to $\mathcal{O} \left(N_e p^{d+1} \right)$.  Hence, the complexity increases by a factor of $p$ going from $d=2$ to $d=3$.  This is one of the major reasons why most spectral element codes use a moderate basis function order $3 \leq p \leq 6$ (see HOMME-NH \cite{taylor2020energy}, NUMA \cite{giraldoetal2013}, etc.).

In 2D, let the semi-infinite domain $\Omega^S$ consist of $N_x$ semi-infinite elements with order $N_{LGL}$ in the finite direction and $N_{LGR}$ in the semi-infinite direction.  
\begin{enumerate}
\item
Cost of constructing a vertical derivative: $\mathcal{O} \left(N_x p_{LGL} p_{LGR}^2 \right)$.
\item
Cost of constructing a horizontal derivative: $\mathcal{O} \left(N_x p_{LGL}^2 p_{LGR} \right)$.  
\end {enumerate}

In 3D, let the semi-infinite domain $\Omega^S$ consist of $N_x N_y$ semi-infinite elements.
\begin{enumerate}
\item
Cost of constructing a vertical derivative: $\mathcal{O} \left(N_x N_y p_{LGL}^2 p_{LGR}^2 \right)$.
\item
Cost of constructing a horizontal derivative: $\mathcal{O} \left(N_x N_y p_{LGL}^3 p_{LGR} \right)$.  
\end {enumerate}

Going from 2D to 3D, the cost of constructing either a vertical or horizontal derivative asymptotically increases by a factor of $N_y p_{LGL}$.  This is the same factor incurred in the finite domain $\Omega^F$ going from 2D to 3D.  Hence, the semi-infinite element approach has a similar scaling with respect to the number of dimensions as the standard tensor product SEM utilizing Lagrange polynomials defined at LGL points. 

\section{Conclusion}\label{sec:conc}
We presented a nodal continuous Galerkin (CG) approach for discretizing PDEs on unbounded domains in one or two dimensions. Utilizing a tensor product Legendre-Laguerre basis, waves/disturbances may propagate between finite and semi-infinite domains without artificial reflections. This coupling is achieved by replacing one of the basis functions in a standard tensor product spectral elements approach with a scaled Laguerre function (SLF) using collocated Legendre-Gauss-Radau (LGR) interpolation/quadrature points. Consequently, a CG spectral element in the finite domain is coupled through the DSS operator to an unbounded domain. In the unbounded domain, a standard Legendre basis on LGL points is used in the finite direction, and a SLF basis with LGR points is used in semi-infinite direction. This approach is employed to construct absorbing layers on the semi-infinite domain $\Omega^S$, showcasing its effectiveness in damping outgoing waves in one dimension with minimal reflection. Subsequently, we demonstrate the ability of this approach to achieve accurate solutions for the advection-diffusion and Helmholtz equations in two dimensions. Finally, we illustrate that this method can effectively build absorbing layers for the non-linear compressible Euler equations in two dimensions. We showcase its capability to damp vertically propagating gravity waves in two standard idealized cases: a linear hydrostatic mountain and the Sch\"ar mountain.  Furthermore, we demonstrate that the semi-infinite element approach is less expensive than the standard Rayleigh damping approach for a variety of 1D and 2D cases, and this efficiency increases with the desired size of the unbounded domain.

Reducing time-to-solution in large high-resolution simulations is crucial. The time typically allocated for computations in the absorbing layer could instead be directed towards the domain of interest or even discarded, favoring quicker simulations. This study does not aim to establish criteria for selecting the appropriate scaling parameter; hence, we limit ourselves to merely altering the domain size for test cases comparable with the literature.  In this regard, the work presented in \cite{Xia2021,Xia2023,Chou2023} opens up the possibility of time-evolving adaptive unbounded domains for simulating atmospheric flows. Additionally, we are extending this model to 3D for efficient moist atmosphere simulations \cite{tissaouiEtAl2023} and expanding this method to cases requiring an absorbing layer in both the lateral and vertical dimensions. 

\section{Code availability} The open source code {\tt Jexpresso} used for this paper is available on Github \cite{jexpressoGithub}. All of the tests described in this paper can be accessed and run from the code at this link:\\
\url{https://doi.org/10.5281/zenodo.11490036}. Timings were obtained by running the code with Julia 1.9.3 on a Macbook Air M1 2020, with macOS Big Sur Version 11.6.
\section{Acknowledgments}
The authors are grateful to the comments by two anonymous reviewers and by Tommaso Benacchio (Danish Meteorological Institute), who all improved the paper.
Tissaoui and Marras gratefully acknowledge the partial support by the National Science Foundation through NSF grant PD-2121367.
Kelly gratefully acknowledges the support of the Office of Naval Research under grant \# N0001419WX00721. 
Makie \cite{DanischKrumbiegel2021} and VisIt \cite{visitWeb} were used to create the plots.

\bibliographystyle{plain}
\bibliography{bibliography_completa}
\newpage
\section*{Appendix A: Constructing the element right hand side}
\label{sct:appendix}
We provide the reader with a pseudo-code for the use of inexact integration in computing the element right hand sides of a PDE on semi-infinite elements, where $\mathbf{rhs}^{Lag}$ is the right hand side a semi-infinite element and $F$ and $G$ are fluxes as defined in (\ref{eq:CL}):
\begin{algorithm}[h!]
    \caption{Construction of the right hand side for an element of the semi-infinite domain}
    \begin{algorithmic}
        \State $\mathbf{rhs}^{Lag} = {\rm zeros}(N_{LGL},N_{LGR})$
        \For {$j =1,N_{LGR}$}
        \For {$i =1,N_{LGL}$}
        \State $\overline{\omega}= \omega(\xi_i)\hat{\omega}(\eta_j)$
        \State $dFd\xi=dFd\eta=dGd\xi=dGd\eta=0$
        \For {$k =1,N_{LGL}$}
        \State $dFd\xi = dFd\xi + h'_k(\xi_i)F(\mathbf{x}(\xi_k,\eta_j,e))$
        \State $dGd\xi = dGd\xi + h'_k(\xi_i)G(\mathbf{x}(\xi_k,\eta_j,e))$
        \EndFor
        \For{$k=1,N_{LGR}$}
        \State $dFd\eta = dFd\eta + \hat{h}^{'Lag}_k(\eta_j)F(\mathbf{x}(\xi_i,\eta_k,e))$
        \State $dGd\eta = dGd\eta + \hat{h}^{'Lag}_k(\eta_j)G(\mathbf{x}(\xi_i,\eta_k,e))$
        \EndFor
        \State $dFdx = dFd\xi \cdot d\xi dx + dFd\eta \cdot d\eta dx$
        \State $dGdz = dGd\xi \cdot d\xi dz + dGd\eta \cdot d\eta dz$
        \State $\mathbf{rhs}^{Lag}_{ij} = \mathbf{rhs}^{Lag}_{ij} - \overline{\omega}|J(\xi_i,\eta_j)|(dFdx + dGdz) $
        \EndFor
        \EndFor
    \end{algorithmic}
\end{algorithm}

With the element right hand sides determined, the same DSS operation described in \S~\ref{sec:num} can be used to construct the global right hand. Similarly to applying DSS to the mass matrix, this enforces the continuity of the global solution and is the only coupling between the finite and semi-infinite domains.

\section*{Appendix B: Extension to 3D}
\label{sct:appendixB}

The proposed semi-infinite element approach may be extended to 3D.  Let $\mathbf{\xi} = (\xi,\eta,\zeta)$ be the coordinate of a point on the three dimensional reference element. Eq. (\ref{eq:basis}) is generalized to three dimensions as follows:

\begin{equation}
    \psi_l(\mathbf{x}) = h_i[\xi(\mathbf{x})] \otimes h_j[\eta(\mathbf{x})] \otimes \overline{h}_k[\eta(\mathbf{x})] , \hspace{5pt} \hspace{3pt} l=i + (j-1)N_{LGL} + (k-1)N_{LGL}^2,
\end{equation}
where $i \in \{1,\dots,N_{LGL}\}$, $j \in \{1,\dots,N_{LGL}\}$, and $k \in \{1,\dots,N_{LGR}\}$.
It is then simple to extend equation (\ref{eq:Gauss_quadrature}) to three dimensions by adding an additional sum over the LGL nodes and including their corresponding weights in the product:

\begin{equation}
    \int_{\Omega_e}f(\mathbf{x})d\mathbf{x} = \int_{\Omega_{ref}}f(\boldsymbol{\xi})|\mathbf{J}(\boldsymbol{\xi})|d\boldsymbol{\xi}\approx \sum_{i,j=1}^{N_{LGL}} \sum_{k=1}^{N_{LGR}}\omega(\xi_i) \omega(\eta_j) \hat{\omega}(\zeta_k)f(\xi_i,\eta_j,\zeta_k)|\mathbf{J}(\xi_i,\eta_j,\zeta_k)|.
\end{equation}

The remainder of the extension, such as the construction of the mass matrix can be done by following the same approach.

\section*{Appendix C: Error estimates for the linear hydrostatic mountain and Sch\"ar mountain cases}

We present here the root mean squared error (RMSE) obtained from comparing two simulations of the linear hydrostatic mountain (LHM). The first is a standard spectral element simulation without any semi-infinite elements and the second uses semi-infinite elements of order 14 in the absorbing layer. The resolutions are the same for the finite domains in both simulations which use spectral elements of order 4. The domain extents are also the same for both simulations. The same comparison is also performed for the Sch\"ar mountain test case but with the order of elements within the finite domain being 10. 

The RMSEs for both cases and all the variables are presented in table \ref{tab:error_LHM_schar} for the solution at $t=10$ hours.

\begin{table}[h!]
    \centering
    \begin{tabular}{|c|c|c|c|c|}
    \hline
      Variable & $\rho$ & $u$ & $v$ & $\theta$ \\
      \hline
        RMSE (LHM) &  2.50e-7 & 3.89e-4 & 3.71e-5& 3.18e-4\\
        \hline
        RMSE (Sch\"ar) & 5.95e-6 & 1.54e-2 & 3.50e-3& 4.80e-3\\
        
        \hline
    \end{tabular}
    \caption{Root mean squared error for each variable at $t=10~$hours. The error is obtained by comparing a simulation using a  standard spectral element approach against a simulation using semi-infinite elements of order 14 within the absorbing layer.}
    \label{tab:error_LHM_schar}
\end{table}

For both mountain wave cases, the order of magnitude of the errors shown in table \ref{tab:error_LHM_schar} are similar or smaller than those obtained by Giraldo and Restelli in 2008 (\cite{giraldo_2008}) who compare their spectral element method against the linearized analytical solutions. This makes us confident that the proposed approach does not introduce any significant errors compared to a standard spectral element approach for these two test cases. 
Furthermore, for LHM case, the total time to solution for the standard spectral element approach was $\sim39,000~$ seconds, while the approach utilizing semi-infinite elements had a total time to solution of $\sim24,000~$ and is thus $\sim1.6$ times faster. For the Sch\"ar mountain case, the total time to solution for the standard spectral element approach was $\sim21,000~$ seconds, while the approach utilizing semi-infinite elements had a time to solution of $\sim13,000~$ seconds and is thus also $\sim1.6$ times faster.



\end{document}